\documentclass{article}
\usepackage{graphicx} % Required for inserting images
\usepackage{subfig}
\usepackage{float}
\usepackage{amsmath} 
\usepackage[url=false,sorting=nyt,sortcites=true,abbreviate=false,giveninits=true]{biblatex}
\AtEveryBibitem{\clearfield{month}}
\AtEveryBibitem{\clearfield{day}}
\usepackage{color}
\usepackage{lipsum}
\usepackage{subfig}
\usepackage{amsfonts}
\usepackage{graphicx}
\usepackage{epstopdf}
\usepackage{algorithmic}
\usepackage{hyperref}
\DeclareGraphicsExtensions{.pdf,.eps,.ps,.jpg,.epsi,.png}

\graphicspath{{images/}{./}}

\newtheorem{Remark}{Remark}

\title{A semi-implicit exactly fully well-balanced relaxation scheme for the Shallow Water Linearized Moment Equations}

\author{C. Caballero-Cárdenas, I. Gómez-Bueno, A. Del Grosso,\\ J. Koellermeier, T. Morales de Luna}

\date{}

\begin{document}

\maketitle

\begin{abstract}
When dealing with shallow water simulations, the velocity profile is often assumed to be constant along the vertical axis. However, since in many applications this is not the case, modeling errors can be significant. Hence, in this work, we deal with the Shallow Water Linearized Moment Equations (SWLME), in which the velocity is no longer constant in the vertical direction, where a polynomial expansion around the mean value is considered. The linearized version implies neglecting the non-linear terms of the basis coefficients in the higher order equations. As a result, the model is always hyperbolic and the stationary solutions can be more easily computed. Then, our objective is to propose an efficient, accurate and robust numerical scheme for the SWLME model, specially adapted for low Froude number situations. Hence, we describe a semi-implicit second order exactly fully well-balanced method. More specifically, a splitting is performed to separate acoustic and material phenomena. The acoustic waves are treated in an implicit manner to gain in efficiency when dealing with subsonic flow regimes, whereas the second order of accuracy is achieved thanks to a polynomial reconstruction and Strang-splitting method. We also exploit a reconstruction operator to achieve the fully well-balanced character of the method.
Extensive numerical tests demonstrate the well-balanced properties and large speed-up compared to traditional methods.

\end{abstract}

\section{Introduction} 
 
Shallow water models have many applications in atmospheric and oceanic dynamics, which have motivated extensive research over the past decades: \cite{VreugdenhilSWE,casulli2000unstructured,castro2001q,gerbeau2001derivation,chinnayya2004well,audusse2004fast,noelle2006well,berthon2012efficient,Chertock2024,XiChen2023,micheldansac2016well,berthon2018simple,busto2022staggered}, among many others. 
These models are commonly represented by the shallow water equations, or the Saint-Venant system in one-dimensional form. The equations can be derived by performing a depth-averaging procedure on the Navier-Stokes equations, which assumes a constant horizontal velocity across the vertical axis. 
While this assumption can simplify calculations in many cases, there are numerous situations where a more complex approach is required. The vertical averaging involved leads to a loss of vertical detail, introducing modelling errors that may be significant and should not be neglected \cite{fernandeznieto2014multilayer,koellermeier2020analysis,garresdiaz2023general}.

Looking at the literature, a first possibility in this regard is to decompose the fluid into two or more layers, each with its own constant horizontal velocity \cite{castro2001q,2layerCastro,bonaventura2018multilayer,fernandeznieto2014multilayer,IFCP2layer,AudusseMultiLayerMassExchange,audusse2011approximation,guerrerofernandez2020second,ADG2layer}. 
%However, this model is derived under the non-physical assumption that the different layers of liquids cannot mix. As a result, the system is not always hyperbolic as the eigenvalues can be complex, leading to instabilities in the numerical simulations. Solutions in this regard have been investigated, by including friction terms \cite{2layerCastro} or terms for the mass exchange \cite{AudusseMultiLayerMassExchange}.

In another direction, a second possibility is provided by shallow water moment models, in which the velocity is no longer constant in the vertical direction as a polynomial expansion around the mean value is considered \cite{kowalski2018moment}. The expansion coefficients, also called moments, allow for the representation of complex vertical velocity profiles. The greater the number of moments, the better the vertical approximation of the velocity.
However,  the model may not be hyperbolic if two (or more) moments are considered. We note that a number of alternative models has been derived in the literature for the 1D case \cite{steffler1993depth,cantero‐chinchilla2018depth,koellermeier2020analysis,garres2021,HuangKoellermeierEQStab} and more recently also for the 2D extension \cite{verbiest2023hyperbolicaxisymmetricshallowwater,bauerle2023rotationalinvariancehyperbolicityshallow}. The technique can also be combined with a multilayer approach by proposing a vertical decomposition into layers that use withing a moment approach \cite{garresdiaz2023general}.

Let us remark that the derivation of steady states for those models is intricate due to the combination of non-conservative products and non-linear terms. 
A solution to this problem is provided in \cite{JulianErnestoLinearizedMoment22}, where the authors propose a partially linearized version of the model which allows for small deviations from the constant velocity profile. The new Shallow Water Linearized Moment Equations (SWLME) are obtained neglecting non-linear terms of the basis coefficients in the higher order equations. Not only this new model is always hyperbolic, regardless of the number of moments, but it also allows for a simpler description of the model's stationary solutions, whereas before they were complicated both analytically and numerically. This fact is of no small importance as it is essential to be able to preserve steady states to deal with small deviations from them, otherwise numerical artifacts and spurious oscillations may appear in the numerical results. 

Numerical preservation of steady states has recently gained a lot of interest. A numerical scheme with this property is defined either \emph{well-balanced} \cite{audusse2004fast,audusse2015simple,diaz2013high,micheldansac2016well,castro2020well} if only zero-velocity steady states are preserved or \emph{fully} well-balanced \cite{LPfullyWB,berthon2016fully,berthon2015fully,berthon2018simple,castro2020well} if all steady states can be preserved. In fact, one may distinguish between well-balanced or exactly well-balanced schemes depending on whether an approximation of the exact steady states are preserved \cite{gomezbueno2021collocation}.

For all the previous reasons, in this paper, we aim to numerically simulate shallow  water flows with non-constant vertical velocity by means of the SWLME model. Then, our objective is to propose an efficient, accurate and robust numerical scheme for the selected model. 
Concerning the last property, we aim for a scheme which is not only consistent, but also exactly fully well-balanced.
Then, the numerical method should be efficient, namely numerical simulations should be fast by neglecting restrictive CFL conditions for the time step. A large part of applications and simulations for shallow water flows occur in subsonic regimes, i.e. when material waves are much slower than acoustic waves.
In this type of situation, the CFL condition may be restrictive due to the fast acoustic waves, which are not of main interest as they quickly leave the domain. 
Hence, we aim to deal with this problem by exploiting a numerical scheme whose CFL condition is based on the slow material waves. This also implies that the transport waves will be better approximated as the CFL condition is adapted to them. For this purpose, we look for an approach that treats implicitly the acoustic waves and explicitly the material waves. Different possibilities arise for this purpose. One option would be to use an IMEX approach as in \cite{bispen2014imex,casulli1994stability,MaccaPhd}. Another alternative is the use of Lagrange-projection methods, which provide a natural decomposition of acoustic and transport waves \cite{StauffertLPimex,chalons2016all,LPerosion,caballerocardenas2023implicit,Bourgeois2024}.
However, even if this strategy was extended to two-layer shallow water flows \cite{ADG2layer}, this was far from trivial and the moment models seem to be even more difficult to deal with in that framework.
A different method for splitting of acoustic and transport phenomena can be found in \cite{Iampietro18}. It was  originally developed for the compressible Euler equations and then extended to the Ripa model \cite{RipaIMEX} and the shallow water system \cite{Celia2024splitting}. It combines a relaxation approach with an splitting of acoustic and transport waves that avoids several difficulties when compared to Lagrange-projection approach.
 
It is important to take into account that implicit schemes can be in general highly diffusive, hence at least second order of accuracy is needed. However, combining implicit methods with a higher order of accuracy could become a complex task, since a naive extension could lead to non-linear systems to be solved iteratively and thus a higher computational cost. Of course, this would be in contradiction with the goal of having an efficient method. One of the advantages of this technique, as shown in \cite{Celia2024splitting}, is precisely the fact that, thanks to the relaxation, only linear systems need to be solved.

In conclusion, we aim to develop a numerical scheme, which is second order accurate and exactly fully well-balanced for the Shallow Water Linearized Moment Equa-tions (SWLME), that splits the acoustic and transport phenomena. The former are then solved implicitly in a efficient way using a relaxation approach. 
%We emphasize that all these complex  issues have already been addressed at least one at a time but, with this work, we aim to combine all of them. 
\\

\emph{Outline of the paper}.
In the next section \ref{sec:Model}, we present the mathematical model and splitting we exploit for the numerical simulations and schemes. Then, both first- and second order methods are described in section \ref{sec:Scheme} along with details for the exactly fully well-balanced character.
Numerical simulations are presented in section \ref{sec:NumSim} to show the good behaviour of the schemes.
Finally, we draw conclusions in section \ref{sec:Conclusion}.

\section{Model description} \label{sec:Model}
The classical shallow water equations (SWE) model the evolution of the water height $h$ and the mean velocity $u_0$ under shallow conditions. They can be derived by depth-averaging of the incompressible Navier-Stokes equations. For the frictionless case with gravitational constant $g$ and bottom topography $z(x)$, the evolution equations are given by
\begin{equation}\label{e:SWE}
    \partial_t
    \begin{pmatrix}
    h\\
    h u_0\\
    \end{pmatrix} +\partial_x
    \begin{pmatrix}
    h u_0\\
    h u_0^2 + \frac{1}{2}g h^2 \\
    \end{pmatrix} =
   \begin{pmatrix}
    0\\
    -gh \partial_x z\\
    \end{pmatrix},
\end{equation}
The simplicity of the model allows to easily compute analytical properties of the model such as its eigenvalues, given by $\lambda_{1,2} = u_0 \pm \sqrt{gh}$, and its steady states given by 
\begin{eqnarray}\label{e:SWE_steady_states}
    h u_0 &=& const,\\
    \frac{1}{2}u_0^2 + g(h+z) &=& const.
\end{eqnarray}

However, the SWE suffer from the inability to model vertical variations of the velocity, since they only include a depth-averaged velocity $u_0$. It has been shown in many applications, that vertical velocity variations are an important feature in shallow flows and need to be taken into account for an accurate representation, see also \cite{kowalski2018moment}. To overcome this deficiency, different extended shallow water models have been derived. The model used in this work is based on the Shallow Water Moment Equations (SWME), which assume an expansion of the velocity profile $u(t,x,\xi)$ around the depth-averaged value:
\begin{equation}\label{e:expansion}
    u(t,x,\xi)=u_0(t,x)+\sum_{i=1}^{N}u_i(t,x)\phi_i\left(\frac{\xi-z(x)}{h(t,x)}\right), \quad \xi \in [z(x),z(x)+h(t,x)]
\end{equation}
where $N$ is referred to as the order of the model, i.e., how many expansion terms are used, $u_i$ are the corresponding expansion coefficients, also called \emph{moments}, and $\phi_i:[0,1]\rightarrow\mathbb{R}$ are \emph{scaled Legendre polynomials} of degree $i$, also called basis functions, defined as
\begin{equation} \label{eq:defphi}
  \phi_i(\zeta) = \frac{1}{i!} \frac{d^i}{d\zeta^i} (\zeta - \zeta^2)^i.
\end{equation}
The expansion \eqref{e:expansion} allows to include polynomial velocity profiles at the expense of including the additional moments $u_i$ as variables. Evolution equations for these moments are derived by computing higher order projections of the Navier Stokes equations \cite{kowalski2018moment}. This results in the SWME, a closed system of $N+2$ evolution equations for the variables $h, u_0, u_1, \ldots, u_N$. Nevertheless, SWME present two main drawbacks: (1) the model is not always hyperbolic, leading to potential instabilities as analysed in detail in \cite{koellermeier2020analysis}, and (2) it is not possible to derive an easy and general expression, using algebraic relations, for steady states as in \eqref{e:SWE_steady_states}, due to the multitude of non-linear terms originating from the expansion and projection.

In this paper, we consider the so-called Shallow Water Linearized Moment Equations (SWLME), derived first in \cite{JulianErnestoLinearizedMoment22} and further analyzed in \cite{HuangKoellermeierEQStab}. The model overcomes both deficiencies of the original SWME model from \cite{kowalski2018moment} using the following assumption: the velocity profile $u(t,x,\xi)$ only has small deviations from a constant profile, i.e., $u_i = \mathcal{O}\left(\epsilon\right)$ for $i\geq 1$. In the higher-order equations, all term of order $\mathcal{O}\left(\epsilon^2\right)$ are then neglected, while the conservation of mass and the momentum equation are left unchanged.  

The SWLME are given by 
\begin{equation}\label{eq:linearised_swe}
    \partial_t \begin{pmatrix}
        h\\
        hu_0\\
        hu_1\\
        \vdots \\
        hu_N
    \end{pmatrix}
    + \partial_x \begin{pmatrix}
        hu_0\\
        hu_0^2+g\frac{h^2}{2}+\frac{1}{3}hu_1^2 + \dots + \frac{1}{2N+1}hu_N^2\\
        2hu_0u_1\\
        \vdots \\
        2hu_0u_N 
    \end{pmatrix}
    = Q \partial_x \begin{pmatrix}
            h\\
        hu_0\\
        hu_1\\
        \vdots \\
        hu_N
        \end{pmatrix} + \begin{pmatrix}
            0 \\
            -gh\partial_x z\\
            0 \\
            \vdots\\
            0
        \end{pmatrix},
\end{equation}
where the non-conservative right-hand side term includes $Q=\text{diag}(0,0,u_0,\dots,u_0).$

As proven in \cite{JulianErnestoLinearizedMoment22}, the eigenvalues are given by 
\begin{equation}\label{SWLME_eigenvalues}
    \lambda_{1,2} = u_0 \pm \sqrt{gh + \sum_{i=1}^N \frac{3 u_i^2}{2i+1}} \quad \textrm{ and } \quad \lambda_{i+2} = u_0, ~\textrm{ for }~ i=1,\ldots,N,
\end{equation}
meaning that the system is hyperbolic for positive water heights. 

The Froude number, defined by  
\begin{equation}\label{e:Froude}
Fr(U)= \displaystyle \frac{|u_0|}{\sqrt{g h}},
\end{equation}
determines the fluid regime: if its value is less than $1$, it is subcritical; if it is $1$ it represents a critical regime, whereas a Froude number bigger than 1 characterises a supercritical flow. 
%In this work we are especially interested in the subcritical regime, since it is in this case in which the improvements in computational efficiency are observed in the implicit schemes with respect to the explicit ones. 
In this work, we are especially interested in the subcritical regime, as this is where the improvements in computational efficiency of implicit versus explicit schemes are observed. 
However, our schemes could be used for other regimes too. 

\begin{Remark}
    The standard definition of the Froude number in equation \eqref{e:Froude} does not take into account the higher moments $u_i$, for $i>0$. In \cite{JulianErnestoLinearizedMoment22}, this has been first discussed by inclusion of additional dimensionless numbers that indicate the deviation from the constant profile. However, since we assume that the expansion coefficients are small, we shall consider here the standard definition of the Froude number from equation \eqref{e:Froude}. Let us remark that one could define an extended Froude number by
    \begin{equation*}
    \widetilde{Fr}(U)=\frac{|u_0|}{ \sqrt{gh + \sum_{i=1}^N \frac{3 u_i^2}{2i+1}}}.
    \end{equation*} 
\end{Remark}

Steady states can be computed by solving
\begin{align}
    h u_0=C_1,\label{C1} \\ 
    \frac{1}{2}u_0^2+g(h+z)+\frac{3}{2}\sum_{i=1}^N \frac{1}{2i+1}u_i^2=C_2, \label{C2}\\ 
    \frac{u_i}{h}=C_k, i\geq 1, k\geq 3, \label{C3}
\end{align}
and we refer to \cite{JulianErnestoLinearizedMoment22} for analysis and discussion on the solvability for existing steady states.

The hyperbolicity and existence of easily computable steady states allows the model to be used in standard numerical simulations. As already said, here we aim at defining a numerical scheme that is well-balanced for the steady states and efficient, specially in the case of low Froude numbers. To this end, we shall consider a relaxation approach for system \eqref{eq:linearised_swe}, given by
\begin{equation}\label{eq:linearised_swe_relaxed}
\begin{aligned}
    \partial_t \begin{pmatrix}
        h\\
        hu_0\\
        hu_1\\
        \vdots \\
        hu_N\\
        h\pi 
    \end{pmatrix}
  &  + \partial_x \begin{pmatrix}
        hu_0\\
        hu_0^2+\pi+\frac{1}{3}hu_1^2 + \dots + \frac{1}{2N+1}hu_N^2\\
        2hu_0u_1\\
        \vdots \\
        2hu_0u_N \\
        a^2 u_0
    \end{pmatrix}
    = \\
    & = 
    \widetilde{Q} \partial_x \begin{pmatrix}
            h\\
        hu_0\\
        hu_1\\
        \vdots \\
        hu_N \\
        h\pi
        \end{pmatrix} + \begin{pmatrix}
            0 \\
            -gh\partial_x z\\
            0 \\
            \vdots\\
            0 \\
            0
        \end{pmatrix} + \begin{pmatrix}
            0\\0\\0\\ \vdots\\0\\ \frac{1}{\varepsilon
} \left( g\frac{h^2}{2} - \pi \right)
        \end{pmatrix},
    \end{aligned}
\end{equation}
where $\widetilde{Q}=\text{diag}(0,0,u_0,\dots,u_0,0)$, the constant $a$ has to fulfill the subcharacteristic condition (see \cite{bouchut2004nonlinear}) given by:
\begin{equation}\label{eq:sub_cond}
 h \sqrt{gh} \leq a.
\end{equation}
Moreover, when $\varepsilon\to 0$, one has $\pi=g\frac{h^2}{2}.$ Here, an instantaneous relaxation of the pressure will be considered as usual. Consequently, the $\varepsilon-$term will be neglected and $\pi$ is then set to an initial pressure value $g\frac{h^2}{2}$ at each time step. Therefore, the resulting model can be written in a compressed way as follows
\begin{equation}\label{eq:pde_compact}
    \partial_t U + \partial_x \left( F(U)+F_P(U)\right) + B(U)\partial_x U + S(U)\partial_x z =0,
\end{equation}
where 
\begin{equation}
    F(U)=(hu_0,hu_0^2+\frac{1}{3}hu_1^2 + \dots + \frac{1}{2N+1}hu_N^2,2hu_0u_1,\dots,2hu_0u_N,0)^T,
\end{equation}
\begin{equation}
    F_P(U)=(0,\pi,0,\dots,0,a^2u_0)^T,
\end{equation}
\begin{equation}
    B=-\widetilde{Q},
\end{equation}
\begin{equation}
    S(U)=(0,gh,0,\dots,0)^T.
\end{equation}

Now, following the ideas in \cite{RipaIMEX,Celia2024splitting}, we perform a splitting, obtaining a pressure system and a transport one.
The pressure system to be solved is written as 
\begin{equation}\label{eq:pressure_system_compact}
    \partial_t U+ \partial_x F_P(U) + S(U)\partial_x z = 0,
\end{equation}
and the one for the transport is given by
\begin{equation}\label{eq:transport_system_compact}
    \partial_t U+ \partial_x F(U)+B(U)\partial_x U=0.
\end{equation}
Note that the eigenvalues of system \eqref{eq:pressure_system_compact} are 
\begin{equation}\label{SWLME_eigenvalues_pressure}
    \lambda_{1,2}^P =  \pm \dfrac{a}{h} \quad \textrm{ and } \quad \lambda_{i+2}^P = 0, ~\textrm{ for }~ i=1,\ldots,N,
\end{equation}
while the ones of system \eqref{eq:transport_system_compact} are 
\begin{equation}\label{SWLME_eigenvalues_transport}
    \lambda_{1,2}^T = u_0 \pm \sqrt{\sum_{i=1}^N \frac{3 u_i^2}{2i+1}} \quad \textrm{ and } \quad \lambda_{i+2}^T
 = u_0, ~\textrm{ for }~ i=1,\ldots,N.
\end{equation}

This means that stability restriction for the pressure system comes from the pressure term, while in the transport system the stability is governed by velocity. Therefore, we shall consider an implicit approach for the pressure system, which will allow to overcome the restrictive CFL condition that comes from the pressure in the case of small Froude number.

\section{Numerical scheme} \label{sec:Scheme}

We will consider a finite volume approach, that is, the spatial domain is divided into a series of cells, denoted by intervals $[x_{i-1/2}, x_{i+1/2})$ for each integer \( i \in \mathbb{Z} \). For simplicity, we assume that each cell has a uniform length, represented by \(\Delta x = x_{i+1/2} - x_{i-1/2}\). Here, the points \( x_{i+1/2} \) mark the interfaces between cells, while the cell centers are designated as \( x_i = \frac{x_{i-1/2} + x_{i+1/2}}{2} \). For the moment, the time variable will remain continuous. Moreover, the approximation of the cell averages will be denoted by
\begin{equation*}
 U_i(t) \approx \frac{1}{\Delta x} \int_{x_{i-1/2}}^{x_{i+1/2}} U(x,t)dx.
\end{equation*}
Since we will consider here schemes up to second order, the cell averages will be identified with the cell center values, which corresponds to a mid point quadrature formula for computing the cell averages. The same quadrature rule will also be used to approximate the integrals of the source term and the non-conservative ones as it will explained later.

In order to obtain fully well-balanced schemes, we will follow the strategy described in \cite{castro2020well}, in which the well-balanced property is transferred to the reconstruction operator. A fully exactly well-balanced reconstruction operator $P_i(x)$ is a reconstruction operator that satisfies 
\[
P_i(x,\{U_j^e\}_{j \in \mathcal{S}_i})=U^e(x), \ \forall x \in [x_{i-1/2},x_{i+1/2}], \ \forall i,
\]
for a given continuous stationary solution $U^e$, and where $\mathcal{S}_i$ denotes the stencil corresponding to cell $i$.

Notice that a stationary solution $U^e$ must satisfy: 
\begin{equation}\label{eq:ode_stationary}
     \partial_x \left( F(U^e)+F_P(U^e)\right) + B(U^e)\partial_x U^e + S(U^e)\partial_x z=0.
\end{equation}
The stationary solutions can be computed by solving \eqref{C1}-\eqref{C3}.
% the following equations:
% \begin{align}
%     hu_0=C_1,\\
%     \frac{1}{2}u_0^2+g(h+z)+\frac{3}{2}\sum_{i=1}^N \frac{1}{2i+1}u_i^2=C_2,\\
%     \frac{u_i}{h}=C_k, i\geq 1, k\geq 3.
% \end{align}
% \textcolor{blue}{ADG: you are repeating the steady states \eqref{C1}-\eqref{C2}-\eqref{C3}.   I would remove the above ones.}

Specifically, being $t_0$ the initial time at each time step, we compute at every time step $[t_0,t_0+\Delta t]$ and at every cell $i$, the steady state $U_i^{e,t_0}$ that satisfies \eqref{eq:ode_stationary} such that $U_i^{e,t_0}(x_i)=U_i^{t_0}$.

For each cell $i$, if the ODE
\eqref{eq:ode_stationary} is subtracted from the system \eqref{eq:pde_compact}, and a similar splitting to the one described above is applied, we can write the semi-discrete version of the pressure and transport systems as follows:
\begin{equation}\label{eq:splitAdetalladorelajado_wb}
\begin{split}
    U_i'(t)=&-\dfrac{1}{\Delta x} \left( F_{P,i+1/2}(t)-F_P(U_i^{e,t_0}(x_{i+1/2}))-F_{P,i-1/2}(t)+F_P(U_i^{e,t_0}(x_{i-1/2})) \right) \\&- \dfrac{1}{\Delta x} \left( \int_{x_{i-1/2}}^{x_{i+1/2}} \left(S(P_{i,U}(x)) - S(U_i^{e,t_0}(x)) \right)z'(x)dx \right),
    \end{split}
\end{equation}
and
\begin{equation}
\begin{split}\label{eq:splitBdetallado_wb}
U_i'(t)=&-\dfrac{1}{\Delta x} \left( F_{i+1/2}(t)-F(U_i^{e,t_0}(x_{i+1/2}))-F_{i-1/2}(t)+F(U_i^{e,t_0}(x_{i-1/2})) \right) \\ 
&-\dfrac{1}{\Delta x} \left( \mathcal{B}_{i+1/2-}(t)+\mathcal{B}_{i-1/2+}(t) \right) \\
&-\dfrac{1}{\Delta x}\left( \int_{x_{i-1/2}}^{x_{i+1/2}} B(P_{i,U}(x,t_0)) \dfrac{dP_{i,U}}{dx} dx-\int_{x_{i-1/2}}^{x_{i+1/2}} B(U_i^{e,t_0}(x)) \dfrac{dU_i^{e,t_0}}{dx} dx \right),
\end{split}
\end{equation}
where 
\begin{itemize}
\item $F_{P,i+1/2}=(0,\pi_{i+1/2}^*,0,\dots,0,a^2u_{0,i+1/2}^*)^T$, being  \(\pi^*_{i \pm \frac{1}{2}}(t)\) and \(u^*_{0,i \pm \frac{1}{2}}(t)\) approximations of the pressure and velocity at the interface, respectively, where \(\pi^*_{i \pm \frac{1}{2}}(t) \approx \pi(x_{i \pm \frac{1}{2}}, t)\) and \(u^*_{0,i \pm \frac{1}{2}}(t) \approx u_0(x_{i \pm \frac{1}{2}}, t)\). 
    \item $P_{i,X}(x,t)$ is the fully well-balanced reconstruction operator of any variable $X$, that we will describe in subsection \ref{subsec:WB_operator}.
    \item $F_{i\pm1/2}=\mathcal{F}(U_{i\pm 1/2-},U_{i\pm 1/2+})$, where $\mathcal{F}$ is a consistent numerical flux, i.e. $\mathcal{F}(U,U)=F(U).$ In this work, the HLL numerical flux is considered (see \cite{harten1997upstream}), that is, 
    \begin{equation}
    \begin{split}
F_{i+1/2}&=\dfrac{F(U_{i+1/2-})+F(U_{i+1/2+})}{2} \\
        &- \dfrac{1}{2} \left(\alpha_{i+1/2}^0 \left( U_{i+1/2+} - U_{i+1/2-} \right) + \alpha_{i+1/2}^1 \left( F(U_{i+1/2+} )- F(U_{i+1/2-}) \right)\right)
    \end{split}
    \end{equation}
    where the coefficients $\alpha_{i+1/2}^0$ and $\alpha_{i+1/2}^1$ are defined as
\begin{equation}
    \alpha_{i+1/2}^0=\frac{S_{i+1/2}^R|S_{i+1/2}^L|-S_{i+1/2}^L|S_{i+1/2}^R|}{S_{i+1/2}^R-S_{i+1/2}^L}, \quad \alpha_{i+1/2}^1=\frac{|S_{i+1/2}^R|-|S_{i+1/2}^L|}{S_{i+1/2}^R-S_{i+1/2}^L}.
\end{equation}
Here, $S_{i+1/2}^L$ and $S_{i+1/2}^R$ are, respectively, approximations of the minimum and maximum wave propagation speeds, computed from equation \eqref{SWLME_eigenvalues}.
\item $\mathcal{B}_{i+1/2\pm}=\dfrac{1}{2}B_{i+1/2}\left( U_{i+1/2+} - U_{i+1/2-} \right) \pm \dfrac{1}{2} \alpha_{i+1/2}^1 B_{i+ 1/2} \left( U_{i+1/2+} - U_{i+1/2-} \right)$, 
where $B_{i+1/2}$ is a consistent matrix with the original matrix $B$. Note that the term $B(U)\partial_xU$ corresponds to a non-conservative product. Therefore, one must give a precise definition of weak solutions. One way to do so is by means of the path-conservative theory (see \cite{pares2006numerical}). Nevertheless, we will not focus on this problem here, which is out of the scope of the paper. In practice, we define
\begin{equation*}
B_{i+1/2} = \int_0^1 B(sU_{i+1/2+}+(1-s)U_{i+1/2-})ds \approx B((U_{i+1/2-}+U_{i+1/2+})/2)
\end{equation*} 
% \begin{equation}
%   B_{i+1/2}\left( U_{i+1/2}^+ - U_{i+1/2}^- \right) \approx \int_0^1 B(\Phi(s;U_{i+1/2}^-,U_{i+1/2}^+))\dfrac{d\Phi}{ds} (s;U_{i+1/2}^-,U_{i+1/2}^+))
% \end{equation}
\end{itemize}

In the previous definitions, the time dependency has been omitted for the sake of simplicity. Moreover, notice that, in system \eqref{eq:splitAdetalladorelajado_wb}, only the variables $hu_0$ and $h\pi$ change with time. Then, using that $h_i$ is constant in this system, we freeze it at each time step when we solve the pressure system. Therefore, for first and second order schemes, taking into account that the mid point rule is applied to approximate the numerical source term, the equations corresponding to the variables $u_0$ and $\pi$ read as follow
\begin{equation}\label{eq:splitAdetalladorelajado_wb_u_pi}
    \begin{cases}
        (u_0)_i'(t)=&-\dfrac{1}{h_i(t_0)\Delta x} \left( \pi_{i+1/2}^*(t)-\pi_{i-1/2}^*(t) - \pi_i^{e,t_0}(x_{i+1/2}) + \pi_i^{e,t_0}(x_{i-1/2}) \right) , \\ 
        \pi_i'(t)=&-\dfrac{a^2}{h_i(t_0)\Delta x} \left( u_{0,i+1/2}^*(t)-u_{0,i-1/2}^*(t)-u_{0,i}^{e,t_0}(x_{i+1/2})+u_{0,i}^{e,t_0}(x_{i-1/2})\right).
    \end{cases}
\end{equation}

Denoting as $\overrightarrow{w}=\pi+au_0$ and $\overleftarrow{w}=\pi-au_0$ the Riemann invariants, system \eqref{eq:splitAdetalladorelajado_wb_u_pi} can be rewritten as two transport equations:
\begin{align}\label{eq:semidisc_pressure_wb_w}
    \begin{cases}
         \overrightarrow{w}_i'(t)&=-\dfrac{a}{h_i(t_0)\Delta x} \left( \overrightarrow{w}_{i+1/2}(t)-\overrightarrow{w}_{i-1/2}(t) - \overrightarrow{w}_i^{e,t_0}(x_{i+1/2}) + \overrightarrow{w}_i^{e,t_0}(x_{i-1/2}) \right)  \\ 
        \overleftarrow{w}_i'(t)&=\dfrac{a}{h_i(t_0)\Delta x} \left( \overleftarrow{w}_{i+1/2}(t)-\overleftarrow{w}_{i-1/2}(t) - \overleftarrow{w}_i^{e,t_0}(x_{i+1/2}) + \overleftarrow{w}_i^{e,t_0}(x_{i-1/2})\right),
    \end{cases}
\end{align}
where \(\overrightarrow{w}_{i+1/2}(t) \approx \overrightarrow{w}(x_{i+1/2}, t)\) and \(\overleftarrow{w}_{i+1/2}(t) \approx \overleftarrow{w}(x_{i+1/2}, t)\) are the numerical fluxes at the intercells, and will be computed by means of a reconstruction operator. In addition, the variables corresponding to the pressure and the velocity can be easily obtained from the following relations:
\begin{equation}\label{eq:rel_w}
    \pi= \dfrac{\overrightarrow{w}+\overleftarrow{w}}{2}, \quad u_0=\dfrac{\overrightarrow{w}-\overleftarrow{w}}{2a}.
\end{equation}
Then, the idea is to first solve \eqref{eq:semidisc_pressure_wb_w} and then  we can use \eqref{eq:rel_w}  to define \(\pi^*_{i+1/2}(t)\) and \(u^*_{0,i+1/2}(t)\) as
\begin{align}
    \pi^*_{i+1/2}(t) &= \dfrac{P_{i,\overrightarrow{w}}(x_{i+1/2}, t) + P_{i+1, \overleftarrow{w}}(x_{i+1/2}, t)}{2}, \label{eq:pi}\\
    u^*_{0,i+1/2}(t) &= \dfrac{P_{i,\overrightarrow{w}}(x_{i+1/2}, t) - P_{i+1, \overleftarrow{w}}(x_{i+1/2}, t)}{2a}. \label{eq:u}
\end{align}

\begin{Remark}
    As far as the stability of the scheme is concerned, for the pressure step the explicit schemes must satisfy the following stability restriction:
\begin{equation}\label{eq:CFL_pressure}
        \Delta t_P \leq \dfrac{\Delta x \max_i \{h_i\}}{a},
    \end{equation}
    while the one for the transport step the restriction can be written as 
\begin{equation}\label{eq:CFL_transport}
     \Delta t_T \leq \dfrac{\Delta x}{2 \max  \{ \max_i \{ |\lambda_{1,i}^T|\},\max_j \{ |\lambda_{2,i}^T|\}   \}}.
    \end{equation}

    Of course, we will choose the maximum $\Delta t$ satisfying both restrictions. The main advantage of using semi-implicit schemes is that they do not require restriction \eqref{eq:CFL_pressure}, which is specially interesting in low Froude number situations.
\end{Remark}

\begin{Remark}
    In practice, the value of the constant $a$ will be computed as 
    \begin{equation*}
        a=\max_i \{h_i \sqrt{gh_i} \}.
    \end{equation*}
\end{Remark}

\subsection{Fully well-balanced reconstruction operator}\label{subsec:WB_operator}

We now illustrate the process for reconstructing the variables, presented in a general form for a variable \(X\), which could represent \(\overrightarrow{w}\), \(\overleftarrow{w}\), \(h\), or \(hu_j\), for $j=0,\dots,N$.

In order to reconstruct a variable $X$ at a given time $t_0$, we apply the strategy described in \cite{castro2020well}, consisting on, at each cell, computing the stationary solution $X_i^{e,t_0}(x)$ and applying a standard reconstruction operator $Q$ to the fluctuations with respect to the stationary solution, giving
\begin{equation}\label{eq:reconstr_t0}
    P_{i,X}^{t_0}(x)=X_i^{e,t_0}(x)+Q_{i,X}^{t_0}(x).
\end{equation}

Following \cite{gomezbueno2023implicit}, the reconstruction operator at time \(t \in [t_0, t_0 + \Delta t]\) is expressed as the sum of the well-balanced reconstruction operator at time \(t_0\) and  a standard reconstruction operator applied to the time fluctuations. Thus, for any variable \(X\), the reconstruction operator is defined as:
\begin{multline*}
P_{i,X}(x,t) = P_{i,X}^{t_0}(x) + \widetilde{Q}_{i,X}(x,t) \\ = X_i^{t_0,e}(x) + Q_{i,X}^{t_0}(x) + \widetilde{Q}_{i,X}(x,t), \quad t \in [t_0, t_0 + \Delta t], \; x \in [x_{i-1/2}, x_{i+1/2}).
\end{multline*}

Here, \(P_{i,X}^{t_0}\) and \(Q_{i,X}^{t_0}\) denote the reconstruction operators defined in \eqref{eq:reconstr_t0} for variable \(X\) at time \(t_0\) and \(\widetilde{Q}_{i,X}(x,t)\) is a reconstruction operator based on time fluctuations, given by 
\[
\widetilde{Q}_{i,X}(x,t) = \widetilde{Q}_i(x; \{X_j^{t,f}\}_{j \in \mathcal{S}_i}), \text{ where } X_j^{t,f} = X_j(t) - X_j^{t_0}, \ j \in \mathcal{S}_i.
\]
%where $\mathcal{S}_i$ denotes the stencil corresponding to cell $i$.
Note that stationary solutions are recalculated at each initial time \(t_0\).

Thus, the first-order well-balanced reconstruction operator is
\begin{equation}\label{reconst_o1}
    P^{o1}_{i,X}(x,t) = X_i^{e,t_0}(x) + X_i(t) - X_i^{e,t_0}(x_i),
\end{equation}
while for second order schemes, it is given by:
\begin{align}\label{reconst_o2}
    P_{i,X}^{o2}(x,t) = X_i^{e,t_0}(x) - X_i^{e,t_0}(x_i) + \Delta X_i^{t_0,f}(x - x_i) + X_i(t) + \Delta X_i^{t,f}(x - x_i).
\end{align}
Here, employing the avg or harmod function limiter from \cite{marquina2007shock}, we set
\begin{equation*}
 \Delta X_i^{t_0,f} = \frac{1}{\Delta x} \left( \phi^{t_0}_{i+} (X_i^{t_0,f} - X_{i-1}^{t_0,f}) + \phi^{t_0}_{i-} (X_{i+1}^{t_0,f} - X_i^{t_0,f}) \right),  
\end{equation*}
where 
\begin{equation*}
    \phi_{i-}^{t_0} = \begin{cases}
             \frac{|d_{i-}|}{|d_{i-}| + |d_{i+}|} & \text{if } |d_{i-}| + |d_{i+}| > 0,\\
             0 & \text{otherwise},
             \end{cases}
\end{equation*}
and 
\begin{equation*}
    \phi_{i+}^{t_0} = \begin{cases}
             \frac{|d_{i+}|}{|d_{i-}| + |d_{i+}|} & \text{if } |d_{i-}| + |d_{i+}| > 0,\\
             0 & \text{otherwise},
             \end{cases}
\end{equation*}
where \(d_{i-} = X_i^{t_0,f} - X_{i-1}^{t_0,f}\) and \(d_{i+} = X_{i+1}^{t_0,f} - X_i^{t_0,f}\), with 
\begin{equation*}
   X_j^{t_0,f} = X_j^{t_0} - X_i^{e,t_0}(x_j)
\end{equation*}
for a given cell \(i\).
Moreover, we define
\begin{equation*}
 \Delta X_i^{t,f} = \frac{1}{\Delta x} \left( \Tilde{\phi}^{t_0}_{i+} (X_i^{t,f} - X_{i-1}^{t,f}) + \Tilde{\phi}^{t_0}_{i-} (X_{i+1}^{t,f} - X_i^{t,f}) \right), 
\end{equation*}
where \(\Tilde{\phi}_{i\pm}^{t_0} = \phi_{i\pm}^{t_0}\) and \(X_i^{t,f} = X_i(t) - X_i^{t_0}\). Note that both \(\Delta X_i^{t_0,f}\) and \(\Delta X_i^{t,f}\) utilize the same limiters calculated at time \(t_0\) to avoid introducing nonlinearities.

\subsection{First order schemes}
We will now describe two different versions of our fully well-balanced first order schemes: one in which both the pressure and the transport systems are solved explicitly, and another one in which the pressure system is solved implicitly, while the transport one is solved explicitly.

Let us assume that the cell averages at time $t_n$, denoted by $U_i^n$, are known. Our goal is to compute the cell averages at the next time step, $U_i^{n+1}$. To accomplish this, we begin by solving the pressure system over a time step $\Delta t$, yielding an intermediate solution, which we denote as $U_i^{n+1-}$. Then, this intermediate solution  $U_i^{n+1-}$ serves as the initial condition for solving the transport system, again over a time step $\Delta t$, which provides the final approximation for $U_i^{n+1}$.

\begin{Remark}
    Note that we could also solve the transport system first, and then the pressure one. However, we have only considered the previously described version of the scheme since in previous works it has been observed that this one performs best in terms of stability (see \cite{Celia2024splitting}).
\end{Remark}

\subsubsection{Explicit scheme}

For the explicit scheme, following the semi-discrete scheme \eqref{eq:splitAdetalladorelajado_wb}, the pressure system updates the variable $hu_0$ as follows:
\begin{equation}\label{eq:exp_o1_pressure}
    (hu_0)_i^{n+1-}=(hu_0)_i^n -\dfrac{\Delta t}{\Delta x} \left( \pi_{i+1/2}^{*,n} - \pi_{i-1/2}^{*,n} - \pi_i^{e,n}(x_{i+1/2}) + \pi_i^{e,n}(x_{i-1/2})\right),
\end{equation}
where $\pi_{i\pm1/2}^{*,n}$ is obtained from the expression \eqref{eq:pi} using the first order reconstruction operator \eqref{reconst_o1} to get the values of the Riemann invariants at the intercells. Note that, within this step $h$ and $(hu_j)$ for $j=1,\dots,N$ do not change, that is, $h_i^{n+1-}=h_i^n$ and $(hu_j)_i^{n+1-}=(hu_j)_i^{n}$ for $j>0$.

Once we have $U_i^{n+1-}$, we can perform the transport step considering \eqref{eq:splitBdetallado_wb} 
\begin{equation}
\begin{split}
U_i^{n+1}=&U_i^{n+1-}-\dfrac{\Delta t}{\Delta x} \left( F_{i+1/2}^{n+1-}-F(U_i^{e,t_{n+1-}}(x_{i+1/2}))-F_{i-1/2}^{n+1-}+F(U_i^{e,t_{n+1-}}(x_{i-1/2})) \right) \\ 
&-\dfrac{\Delta t}{\Delta x} \left( \mathcal{B}_{i+1/2-}^{n+1-}+\mathcal{B}_{i-1/2+}^{n+1-} \right) \\
&-\dfrac{\Delta t}{\Delta x} B(U_i^{n+1-}) \left(U_{i+1/2-}^{n+1-} - U_{i}^{e,n+1-}(x_{i+1/2} ) -U_{i-1/2+}^{n+1-} + U_{i}^{e,n+1-}(x_{i-1/2})\right) ,
\end{split}
\end{equation}
where the reconstructed states are always computed using the well-balanced first order reconstruction operator \eqref{reconst_o1}, the terms $F_{i+1/2}^{n+1-}$ and $\mathcal{B}_{i+1/2-}^{n+1-}$ are evaluated at the reconstructed states at time $t_{n+1-}$ and the derivative in the integral has been approximated using 
\begin{equation*}
    \dfrac{dP_{i,U}}{dx}(x_i) \approx \dfrac{P^{o1}_{i,U}(x_{i+1/2})-P^{o1}_{i,U}(x_{i-1/2})}{\Delta x} = \dfrac{U_{i+1/2-}-U_{i-1/2+}}{\Delta x}.
\end{equation*}

\subsubsection{Semi-implicit scheme}

In the semi-implicit case, the only difference with the explicit one happens in the pressure step. Here, the variable $hu_0$ is updated using 
\begin{equation}\label{eq:imp_o1_pressure}
    (hu_0)_i^{n+1-}=(hu_0)_i^n -\dfrac{\Delta t}{\Delta x} \left( \pi_{i+1/2}^{*,n+1-} - \pi_{i-1/2}^{*,n+1-} - \pi_i^{e,n}(x_{i+1/2}) + \pi_i^{e,n}(x_{i-1/2})\right).
\end{equation}
In practice, this value is obtained as 
\begin{equation*}
    (hu_0)_i^{n+1-}= h_i^n \cdot \dfrac{\overrightarrow{w}_i^{n+1-}-\overleftarrow{w}_i^{n+1-}}{2a},
\end{equation*}
where the variables $\overrightarrow{w}_i^{n+1-}$ and $\overleftarrow{w}_i^{n+1-}$ are given by system \eqref{eq:semidisc_pressure_wb_w} which is discretised as follows
\begin{equation*}
\overrightarrow{w}_i^{n+1-}=\overrightarrow{w}_i^n - \dfrac{a\Delta t}{h_i^n\Delta x} \left(\overrightarrow{w}_{i+1/2}^{n+1} - \overrightarrow{w}_{i-1/2}^{n+1} - \overrightarrow{w}_{i}^{e,n}(x_{i+1/2}) + \overrightarrow{w}_{i}^{e,n}(x_{i-1/2}) \right),
\end{equation*}
and
\begin{equation*}
\overleftarrow{w}_i^{n+1-}=\overleftarrow{w}_i^n + \dfrac{a\Delta t}{h_i^n\Delta x} \left(\overleftarrow{w}_{i+1/2}^{n+1} - \overleftarrow{w}_{i-1/2}^{n+1} - \overrightarrow{w}_{i}^{e,n}(x_{i+1/2}) + \overrightarrow{w}_{i}^{e,n}(x_{i-1/2}) \right),
\end{equation*}
where \eqref{reconst_o1} is used to compute
$\overrightarrow{w}_{i\pm1/2}^{n+1}$  and 
\begin{equation*}
\overrightarrow{w}_{i}^{e,n}(x_{i\pm 1/2})=\dfrac{1}{2}g(h_i^e)^2(x_{i\pm 1/2})+au_{0,i}^e(x_{i\pm 1/2}).
\end{equation*}

\subsection{Second order schemes}

To achieve second order accuracy, we employ the Strang splitting method (see \cite{strang1968construction, marchuk1988metody, marchuk1990splitting}) for time-stepping. Specifically, we proceed as follows:

\begin{enumerate}
    \item Take a half-step for the transport system using a time step of $\Delta t/2$.
    
    \item Advance the pressure system over a full time step $\Delta t$.

    \item Conclude with another half-step for the transport system, again using a time step of $\Delta t/2$, to obtain the final values ${U}_i^{n+1}$.
\end{enumerate}
This could also be written in terms of the approximate solution operators of the systems as
 \begin{equation}\label{eq:strang_TPT}
    U(x,t+\Delta t)= S_T^{\frac{\Delta t}{2}} \circ S_P^{\Delta t} \circ S_T^{\frac{\Delta t}{2}}(U(x,t)),
\end{equation}
where the terms $S_P^{\tau}$ and $S_T^{\tau}$  denote the approximate solution operators over the interval $[t, t+\tau]$ for the pressure system and the transport one, respectively. 

In each step, second order spatial approximations are required. While each sub-step is first order in time, second order accuracy overall is achieved through the Strang splitting method.

\begin{Remark}
    As mentioned for the first order schemes, in this second order case, we could also consider another variant of the scheme:
     \begin{equation*}\label{eq:strang_PTP}
    U(x,t+\Delta t)= S_P^{\frac{\Delta t}{2}} \circ S_T^{\Delta t} \circ S_P^{\frac{\Delta t}{2}}(U(x,t)).
\end{equation*}
However, following again the work \cite{Celia2024splitting}, we have only considered the approach \eqref{eq:strang_TPT} due to stability reasons. 
\end{Remark}

For the sake of readability, the explicit description of the discrete version of the second order schemes are not detailed here. The procedure in this case is analogous to that already described in the first order case using \eqref{eq:strang_TPT} as well as second order reconstructions given in \eqref{reconst_o2}.

\section{Numerical experiments} \label{sec:NumSim}

In this section, we present several numerical experiments designed to evaluate the performance of our proposed schemes. In each experiment we consider open boundary conditions and gravitational acceleration fixed at $g = 9.812$. Unless stated otherwise, the number of moments is set to $N = 8$.

In the following tests, when we mention a CFL number greater than 1, we are specifically referring to the stability constraint related to the sound speed \eqref{eq:CFL_pressure}. Let us remark that our semi-implicit methods also present a stability constraint linked to the transport step, which depends on the flow velocity \eqref{eq:CFL_transport}. This means that we cannot increase the sound-speed-related CFL indefinitely, as it remains bounded by the transport CFL limit.

\subsection{Test 1. Well-balanced property check: lake-at-rest}
In order to check the well-balanced property, let us consider a lake-at-rest initial condition given in \cite{JulianErnestoLinearizedMoment22}. The bottom topography is given by
\begin{equation}\label{test1_bottom}
    z(x)=\left\{ \begin{array}{ll} 2-x^2 & \text{if }  -0.5 \leq x \leq 0.5, \\  
    1.75 & \text{otherwise.} \end{array} \right.
\end{equation}
The water depth is $h(x)=3-z(x)$ and a zero velocity profile is considered 
$u_j(x)=0,$ for $j=0,1,\ldots, 8$.

The spatial domain is $[-1, 1]$ and the experiment is run until the final time $t=0.5$, considering a $400-$cell mesh. We take the CFL value to be $0.9$ for the explicit schemes whereas the implicit methods use the CFL parameter set to $10$.
% We consider the zero velocity profile corresponding to water-at-rest given by the choices:
%\begin{equation}\label{test1_cini}
%    U^0(x)= \left(3-z(x), 0, \dots, 0 \right),
%\end{equation}
%which is shown in Figure \ref{fig:test1_ci}. 
Table \ref{test1_tablaerrores} shows the $L^1-$errors at  $t = 0.5$ between the initial water-at-rest state and the final numerical solutions for the different schemes. For simplicity, in the previous table, only the variables $h$ and $hu_0$ are shown: the errors for the other variables are exactly zero for this test. As expected, all methods give machine accuracy.

\begin{table}[H]
\centering
\begin{tabular}{|c|c|c|c|c|}
\hline
Variable & EXP O1 & EXP O2 & IMP O1 & IMP O2  \\ \hline
$h$ & $1.33\cdot 10^{-15}$ & $4.44\cdot 10^{-16}$ & $3.10\cdot 10^{-15}$&  $2.22\cdot 10^{-15}$\\ \hline
$hu_0$ & $7.48\cdot 10^{-15}$  & $2.06\cdot 10^{-14}$ & $7.60\cdot 10^{-15}$ & $7.60\cdot 10^{-15}$ \\\hline
\end{tabular}
\caption{Test 1. Well-balanced property check: lake-at-rest. $L^1-$errors at time $t = 0.5$. The errors of all schemes are machine precision, indicating the well-balanced property.}
%for the initial condition \eqref{test1_cini}.}
\label{test1_tablaerrores}
\end{table}

%\begin{figure}
%    \centering
%\includegraphics[width=0.75\linewidth]{images/test1_ci.pdf}
%    \caption{Test 1. Well-balanced property check: lake-at-rest. Initial condition}
%    \label{fig:test1_ci}
%\end{figure}

\subsection{Test 2. Well-balanced property check: subcritical stationary flow with zero moments}

\begin{figure}
    \centering
     \subfloat[Free surface]{
\includegraphics[width=0.5\linewidth]{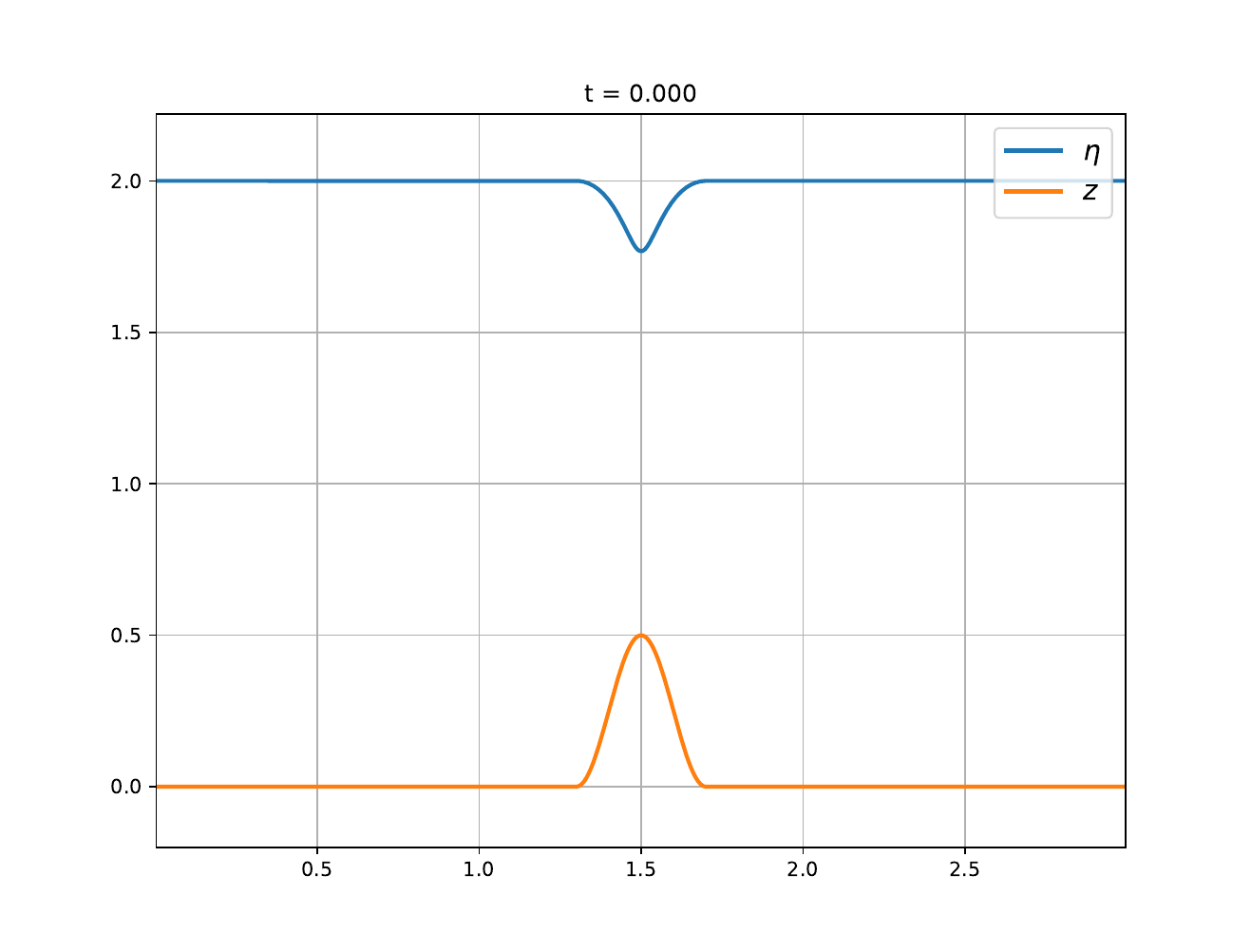}}
\subfloat[Velocity $u_0$]{
\includegraphics[width=0.5\linewidth]{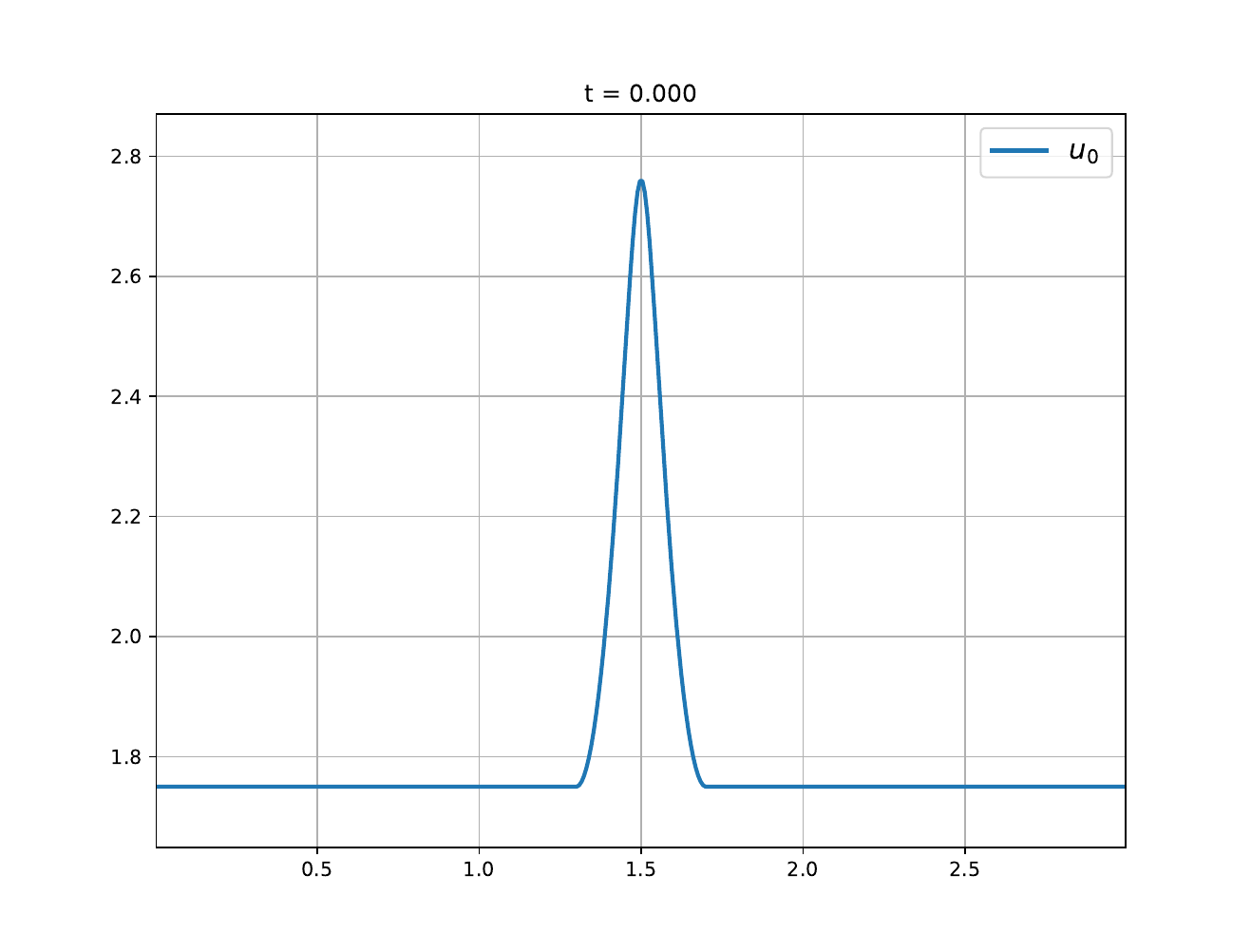}}
    \caption{Test 2. Well-balanced property check: subcritical stationary flow with zero moments. Initial condition: free surface and bottom (left) and velocity $u_0$ (right).}
    \label{fig:test2_ci}
\end{figure}

Let us check the fully well-balanced property, considering now an initial condition which corresponds to a subcritical stationary solution given in \cite{diaz2013high}. The bottom topography is considered to be
\begin{equation}\label{test2_bottom}
    z(x)=\left\{ \begin{array}{ll} \frac{1}{4}\left( \cos{\left( (x+\frac{1}{2}) 5 \pi\right) +1}\right) & \text{if }  1.3 \leq x \leq 1.7, \\  
    0 & \text{otherwise.} \end{array} \right.
\end{equation}
We consider the subcritical steady state characterised by the constants
\begin{equation}\label{test2_cini}
    C_1=3.5, \quad C_2=21.15525, \quad C_k=0, \, k \geq 3,
\end{equation}
in \eqref{C1}-\eqref{C2}-\eqref{C3}, which is shown in Figure \ref{fig:test2_ci}.
The spatial domain is $[0, 3]$ and the experiment is run until  $t=0.5$, considering a $400-$cell mesh. For the explicit and the implicit schemes, we take  the CFL value to be $0.9$ and $1.26$, respectively. Recall that, in practice, the CFL number is restricted by the transport step of the scheme in the implicit case. Since the Froude number varies between $0.4$ and $0.78$, no larger CFL condition is possible in this test.

Table \ref{test2_tablaerrores} shows the $L^1-$errors at time $t = 0.5$ between the initial subcritical initial condition and the final numerical solutions for the different schemes. As expected, all methods give machine accuracy.

\begin{table}[H]
\centering
\begin{tabular}{|c|c|c|c|c|}
\hline
Variable & EXP O1 & EXP O2 & IMP O1 & IMP O2  \\ \hline
$h$ & $2.90\cdot 10^{-14}$ & $3.41\cdot 10^{-14}$ & $3.66\cdot 10^{-14}$&  $3.24\cdot 10^{-14}$\\ \hline
$hu_0$ & $6.26\cdot 10^{-14}$  & $1.82\cdot 10^{-14}$ & $2.35\cdot 10^{-14}$ & $2.08\cdot 10^{-14}$ \\\hline
\end{tabular}
\caption{Test 2. Well-balanced property check: subcritical stationary flow with zero moments. $L^1-$errors at time $t = 0.5$ for the subcritical
initial condition given by the constants \eqref{test2_cini}. The errors of all schemes are machine precision, indicating the well-balanced property.}
\label{test2_tablaerrores}
\end{table}

As mentioned before, the numerical scheme described here is specially designed for situations where the Froude number is low, as it then makes sense to make implicit the pressure part of the system in order to avoid an overly restrictive stability condition. Let us now consider a subcritical steady state scenario with a lower Froude number, characterized by the following constants:
\begin{equation}\label{test2_cini2}
    C_1=0.5, \quad C_2=21.15525, \quad C_k=0, \, k \geq 3,
\end{equation}
Again, the final time is $t = 0.5$. In this case the Froude number is between $0.05$ and $0.07$, which allows us to consider CFL up to $10$ for the implicit scheme. This initial condition is preserved up to machine precision as well, but we are now interested in the efficiency of the scheme. Table \ref{tab:test2_cputime} shows the different CPU times needed for the explicit and implicit schemes. Notice that the first order implicit scheme shows a speedup of $8.45$ compared to the explicit one, while for the second order case the speedup of the implicit scheme is $10.93$ compared to the explicit one. These numbers are close to the relative increase in the CFL number of $11.1$.

\begin{table}
    \centering
    \begin{tabular}{|c|c|c|c|}
    \hline
       EXP O1  & EXP O2  & IMP O1  &  IMP O2\\ \hline
        23.58 & 90.51 & 2.79 & 8.28\\ \hline
    \end{tabular}
    \caption{Test 2. Subcritical stationary flow with zero moments. CPU runtime in s for the different schemes using a $400-$cell mesh. The CFL is equal to $0.9$ for the explicit schemes and $10$ for the implicit ones. As expected by the relative CFL numbers, we observe a speedup of about 10 for the respective implicit schemes.}
    \label{tab:test2_cputime}
\end{table}

\subsection{Test 3. Well-balanced property check: subcritical stationary flow with non-zero moments}

\begin{figure}[H]
 \begin{center}
     \subfloat[ Free surface $\eta$]{
\includegraphics[width=0.5\textwidth]{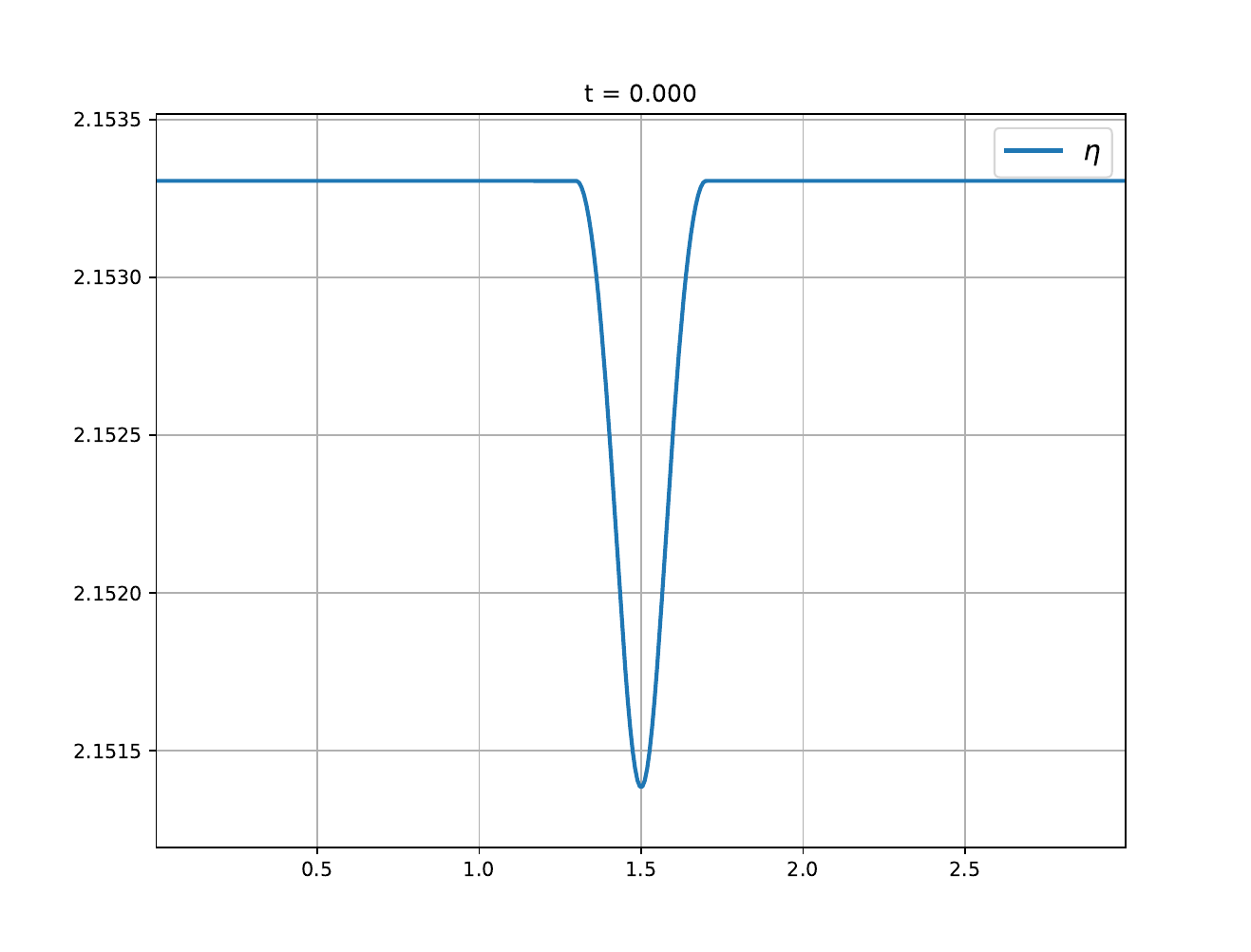}}
\subfloat[ $u_0$ ]{
\includegraphics[width=0.5\textwidth]{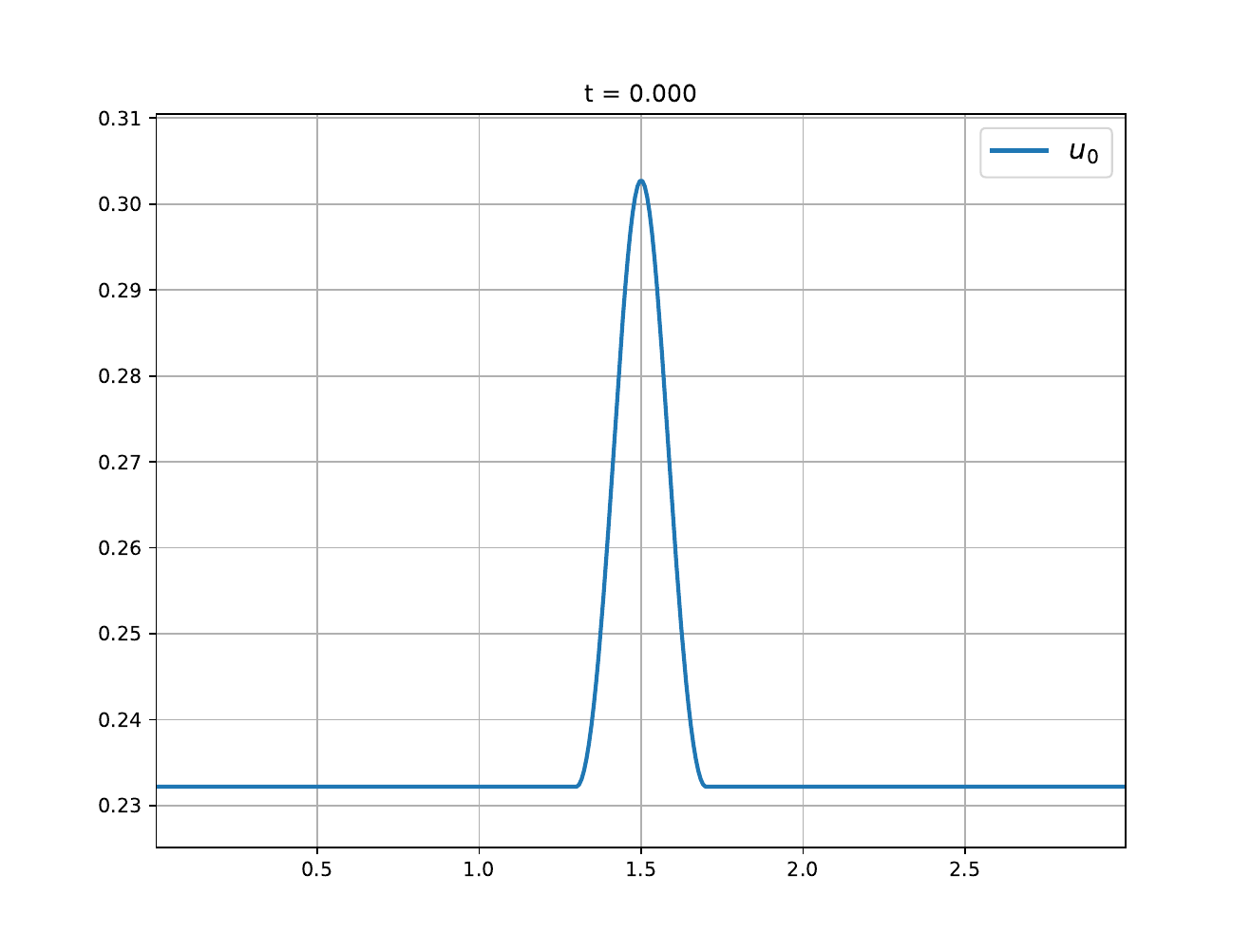}}
\vspace{2mm}
\subfloat[ $u_j, j\geq 1$ ]{
\includegraphics[width=0.5\textwidth]{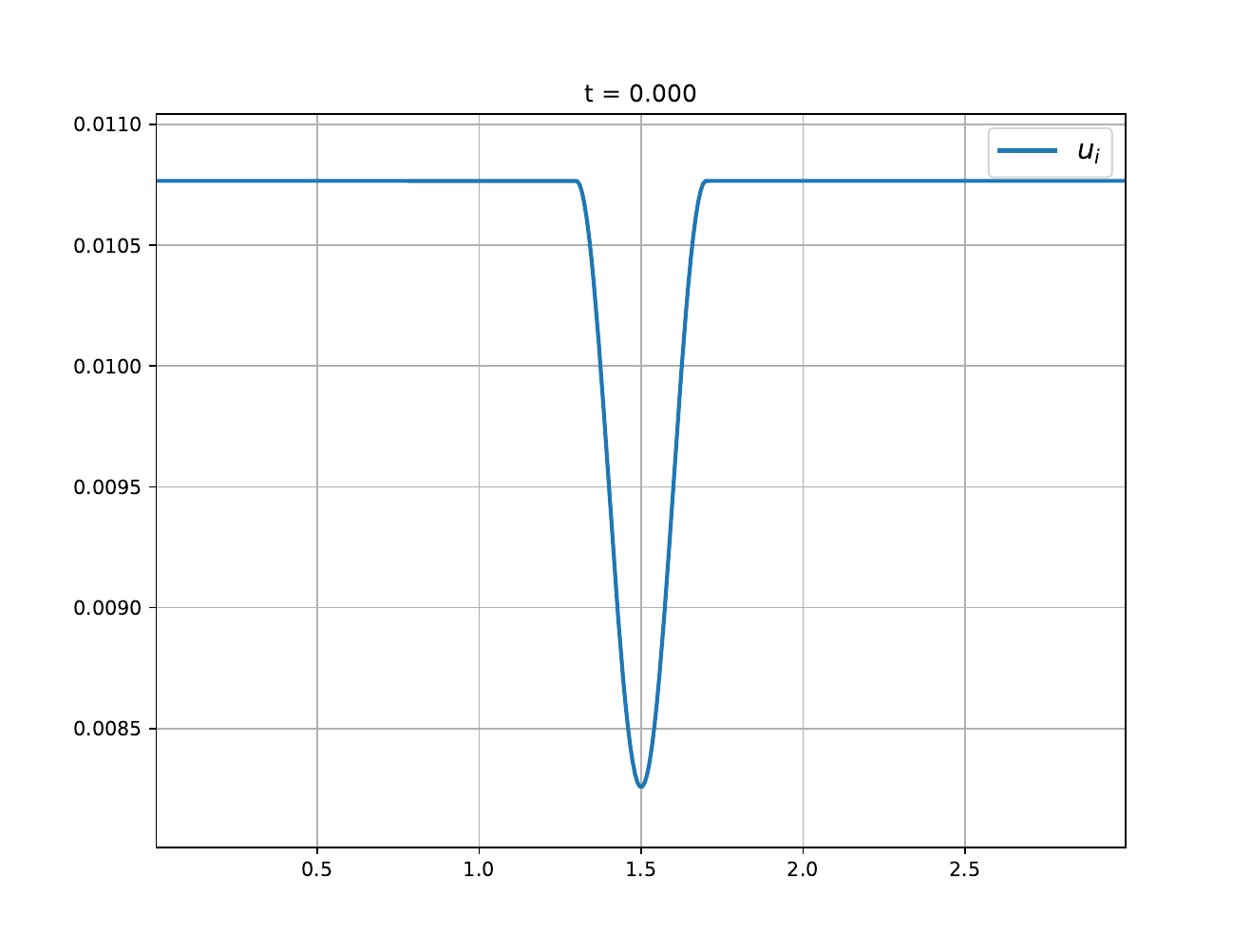}}
\end{center}
    \caption{Test 3. Well-balanced property check: subcritical steady state with non-null moments. Initial condition: free surface $\eta$ (up-left), $u_0$ (up-right) and $u_j, \, j \geq 1$ (down).}
    \label{fig:test3_ci}
\end{figure}

\begin{figure}[H]
    \centering
\includegraphics[width=0.75\linewidth]{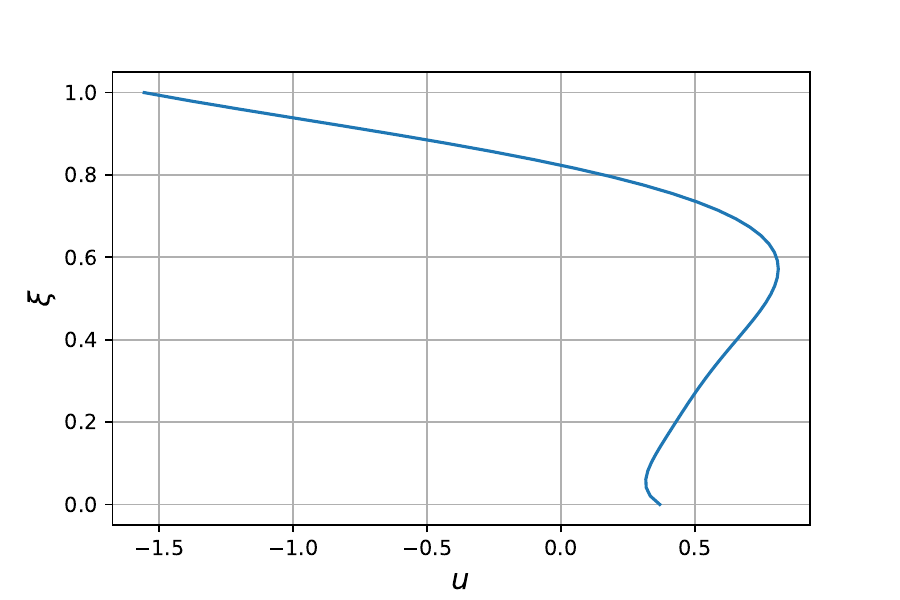}
    \caption{Test 3. Subcritical stationary flow with non-zero moments: velocity profile at $x=1.5$.}  
    %\textcolor{blue}{(JK: it seems this profile has recirculation due to negative u at the top. Might be critizised as non-physical, maybe can be changed by just changing one of the constants $C_i$ above, like $C_3$.)} \textcolor{red}{CCC and IGB: we have considered this test since it is the one used in your paper with Ernesto. We think that in order to check our numerical schemes it is important to compare with previous results existing in the literature. \JK: ok, that is good then!}} 
    \label{fig:test3_velprof}
\end{figure}

Let us consider an initial condition which corresponds to a  stationary solution in the subcritical regime with non-zero moments inspired by an experiment given in \cite{JulianErnestoLinearizedMoment22}. The bottom topography is given again by \eqref{test2_bottom},
the spatial domain is $[0, 3]$ and the final time is $t=0.5$, considering a $400-$cell mesh. We take the CFL value to be $0.9$ for the explicit schemes and $9.15$ for the implicit ones (which is the maximum allowed for this particular solution). We consider the subcritical steady state with non-zero moments characterised by the constants:
\begin{equation}\label{test3_cini}
    C_1=0.5, \quad C_2=21.15525, \quad C_k=0.005, \, k \geq 3,
\end{equation}
which is shown in Figure \ref{fig:test3_ci}. The Froude number is small, between $0.05$ and $0.07$. Figure \ref{fig:test3_velprof} shows the velocity profile at $x=1.5$. Table \ref{test3_tablaerrores} shows the $L^1-$errors at  $t = 0.5$ between the initial subcritical state and the final numerical solutions for the different schemes. As expected, all methods give machine accuracy.   The CPU times for the different schemes are presented in Table \ref{tab:test3_cputime} with the purpose of analysing the efficiency of the methods. Notice that the first and second order implicit schemes show a speedup of $9.4$ and $9.5$, respectively, over the explicit ones. This corresponds directly to the increase of CFL number by a factor of about $10$.

\begin{table}[H]
\centering
\begin{tabular}{|c|c|c|c|c|}
\hline
Variable & EXP O1 & EXP O2 & IMP O1 & IMP O2  \\ \hline
$h$ & $2.66\cdot 10^{-14}$ & $1.68\cdot 10^{-14}$ & $3.46\cdot 10^{-14}$&  $2.02\cdot 10^{-14}$\\ \hline
$q_0$ & $8.03\cdot 10^{-14}$  & $3.73\cdot 10^{-14}$ & $1.46\cdot 10^{-13}$ & $2.99\cdot 10^{-14}$ \\\hline
$hu_j, \, j \geq 1$ & $4.78\cdot 10^{-16}$  & $3.36\cdot 10^{-16}$ & $7.49\cdot 10^{-16}$ & $3.22\cdot 10^{-16}$ \\\hline
\end{tabular}
\caption{Test 3. Well-balanced property check: subcritical stationary flow with non-zero moments. $L^1-$errors at time $t = 0.5$ for the subcritical
initial condition given by the constants \eqref{test3_cini}. The errors of all schemes are machine precision, indicating the well-balanced property.}
\label{test3_tablaerrores}
\end{table}

\begin{table}
    \centering
    \begin{tabular}{|c|c|c|c|}
    \hline
       EXP O1  & EXP O2  & IMP O1  &  IMP O2\\ \hline
        23.07 & 74.85 & 2.46 & 7.86\\ \hline
    \end{tabular}
    \caption{Test 3. Subcritical stationary flow with non-zero moments. CPU times for the different schemes using a $400-$cell mesh. The CFL is equal to $0.9$ for the explicit schemes and $9.15$ for the implicit ones. Again, the speedup is about $10$, in line with the increase of the CFL number.}
    \label{tab:test3_cputime}
\end{table}

\subsection{Test 4. Second order accuracy check}
This experiment is devoted to check the accuracy of the explicit and implicit second order schemes. 
Let us consider a perturbation of the subcritical stationary solution $U^{0,*}$ considered in Test 3, characterised by the constants \eqref{test3_cini}, with the bottom \eqref{test2_bottom}. The initial perturbation is given by
\begin{equation}\label{test4_ci}
    U^0(x)=U^{0,*}(x)+\left(10^{-4}e^{-200(x-2)^2},0, \dots, 0 \right).
\end{equation}
The spatial domain is $[0,3]$ and the final time is $t = 0.1$. We consider four grids consisting of 25, 50, 100 and 200 uniform points in order to compute the errors, comparing with a reference solution, which is computed with a mesh of 800 cells using the same numerical scheme. The CFL number is set to $0.9$ for explicit schemes, whereas a CFL value equal to $2$ is considered for implicit ones. The results are shown in Tables \ref{test4_tablaorden_exp} and \ref{test4_tablaorden_imp_cfl2}:  the expected second order of accuracy is achieved in all variables.

\begin{table}[H]
\centering
\begin{tabular}{|c|cc|cc|cc|}
\hline
No. of & \multicolumn{2}{c|}{$h$}           & \multicolumn{2}{c|}{$u_0$}         & \multicolumn{2}{c|}{$u_j, j\geq 1$} \\ \cline{2-7} 
cells  & \multicolumn{1}{c|}{Error} & Order & \multicolumn{1}{c|}{Error} & Order & \multicolumn{1}{c|}{Error}  & Order \\ \hline
  25 & 7.08e-03 & - & 1.36e-03 & - & 3.54e-05 & - \\ \hline
    50 & 2.09e-03 & 1.76 & 3.07e-04 & 2.14 & 1.05e-05 & 1.76 \\ \hline
   100 & 6.08e-04 & 1.78 & 9.14e-05 & 1.75 & 3.05e-06 & 1.78 \\ \hline
   200 & 1.37e-04 & 2.15 & 2.17e-05 & 2.07 & 6.87e-07 & 2.15 \\ \hline
\end{tabular}
\caption{Test 4. Second order accuracy check. Explicit scheme. CFL$=0.9$.}
\label{test4_tablaorden_exp}
\end{table}

% \begin{table}[H]
% \centering
% \begin{tabular}{|c|cc|cc|cc|}
% \hline
% No. of & \multicolumn{2}{c|}{$h$}           & \multicolumn{2}{c|}{$u_0$}         & \multicolumn{2}{c|}{$u_i, i\geq 1$} \\ \cline{2-7} 
% cells  & \multicolumn{1}{c|}{Error} & Order & \multicolumn{1}{c|}{Error} & Order & \multicolumn{1}{c|}{Error}  & Order \\ \hline
%  25 & 7.09e-03 & - & 1.36e-03 & - & 3.54e-05 & - \\ \hline
% %    25 & 7.08e-03 & 0.00 & 1.36e-03 & 0.00 & 3.54e-05 & 0.00 \\ \hline
%     50 & 2.10e-03 & 1.75 & 3.18e-04 & 2.10 & 1.05e-05 & 1.75 \\ \hline
%    100 & 6.12e-04 & 1.78 & 9.80e-05 & 1.70 & 3.07e-06 & 1.78 \\ \hline
%    200 & 1.39e-04 & 2.14 & 2.51e-05 & 1.96 & 6.97e-07 & 2.14 \\ \hline
% \end{tabular}
% \caption{Test 4. Second order accuracy check. Implicit scheme. CFL=$1$.}
% \label{test4_tablaorden_imp_cfl1}
% \end{table}

% \begin{table}[H]
% \centering
% \begin{tabular}{|c|cc|cc|cc|}
% \hline
% No. of & \multicolumn{2}{c|}{$h$}           & \multicolumn{2}{c|}{$u_0$}         & \multicolumn{2}{c|}{$u_i, i\geq 1$} \\ \cline{2-7} 
% cells  & \multicolumn{1}{c|}{Error} & Order & \multicolumn{1}{c|}{Error} & Order & \multicolumn{1}{c|}{Error}  & Order \\ \hline
%   25 & 7.08e-03 & - & 1.36e-03 & - & 3.54e-05 & - \\ \hline
%     50 & 2.10e-03 & 1.75 & 3.17e-04 & 2.10 & 1.05e-05 & 1.75 \\ \hline
%    100 & 6.12e-04 & 1.78 & 9.93e-05 & 1.67 & 3.07e-06 & 1.78 \\ \hline
%    200 & 1.39e-04 & 2.13 & 2.67e-05 & 1.89 & 7.02e-07 & 2.13 \\ \hline
% \end{tabular}
% \caption{Test 4. Second order accuracy check. Implicit scheme. CFL=$1.5$.}
% \label{test4_tablaorden_imp_cfl1.5}
% \end{table}

\begin{table}[H]
\centering
\begin{tabular}{|c|cc|cc|cc|}
\hline
No. of & \multicolumn{2}{c|}{$h$}           & \multicolumn{2}{c|}{$u_0$}         & \multicolumn{2}{c|}{$u_j, j\geq 1$} \\ \cline{2-7} 
cells  & \multicolumn{1}{c|}{Error} & Order & \multicolumn{1}{c|}{Error} & Order & \multicolumn{1}{c|}{Error}  & Order \\ \hline
  25 & 7.08e-03 & - & 1.36e-03 & - & 3.54e-05 & - \\ \hline
    50 & 2.10e-03 & 1.75 & 3.16e-04 & 2.10 & 1.05e-05 & 1.75 \\ \hline
   100 & 6.14e-04 & 1.77 & 1.04e-04 & 1.61 & 3.08e-06 & 1.77 \\ \hline
   200 & 1.41e-04 & 2.13 & 2.93e-05 & 1.82 & 7.08e-07 & 2.12 \\  \hline
\end{tabular}
\caption{Test 4. Second order accuracy check. Implicit scheme. CFL=$2$.}
\label{test4_tablaorden_imp_cfl2}
\end{table}

% \begin{figure}[H]
%  \begin{center}
%      \subfloat[ Explicit scheme]{
% \includegraphics[width=0.49\textwidth]{images/orden2_explicito.pdf}}
% \subfloat[ Implicit scheme. CFL$=1$.]{
% \includegraphics[width=0.49\textwidth]{images/orden2_implicito_CFL1.pdf}}
% \vspace{2mm}
% \subfloat[ Implicit scheme. CFL$=1.5$.]{
% \includegraphics[width=0.49\textwidth]{images/orden2_implicito_CFL1.5.pdf}}
% \subfloat[ Implicit scheme. CFL$=2$.]{
% \includegraphics[width=0.49\textwidth]{images/orden2_implicito_CFL2.pdf}}
% \end{center}
%     \caption{Test 4. Second order accuracy check. Explicit scheme (up-left) and implicit scheme with CFL$=1$ (up-right), CFL$=1.5$ (down-left) and CFL$=2$ (down-right).  \textcolor{blue}{(JK: I think either axis labels or description of axes in the caption would be necessary.)}}
%     \label{fig:test4_orden}
% \end{figure}

\subsection{Test 5. Perturbation of a subcritical steady state}

\begin{figure}[H]
 \begin{center}
     \subfloat[ Free surface $\eta$]{
\includegraphics[width=0.5\textwidth]{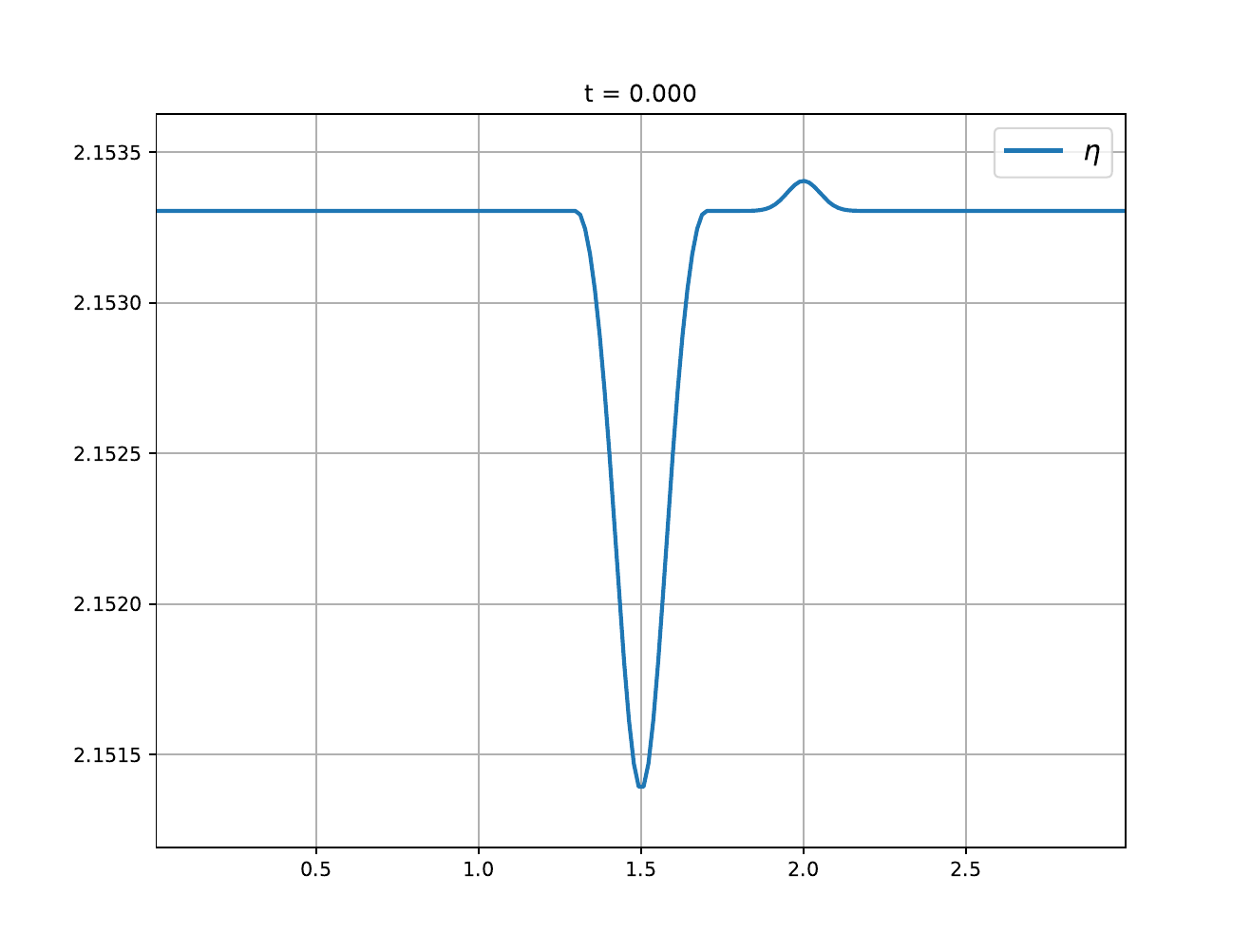}}
\subfloat[ $hu_0$ ]{
\includegraphics[width=0.5\textwidth]{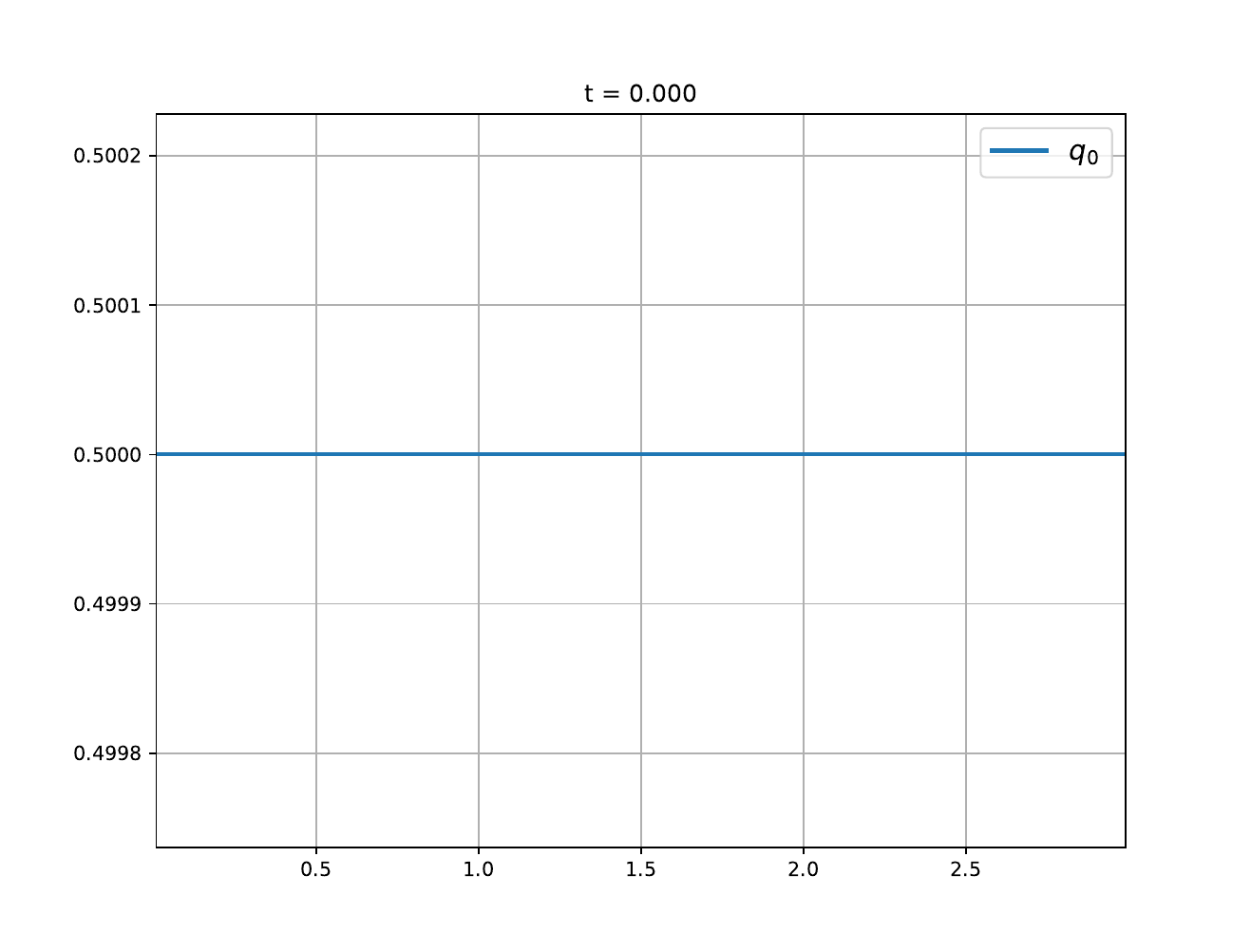}}
\vspace{2mm}
\subfloat[ $hu_1$ ]{
\includegraphics[width=0.5\textwidth]{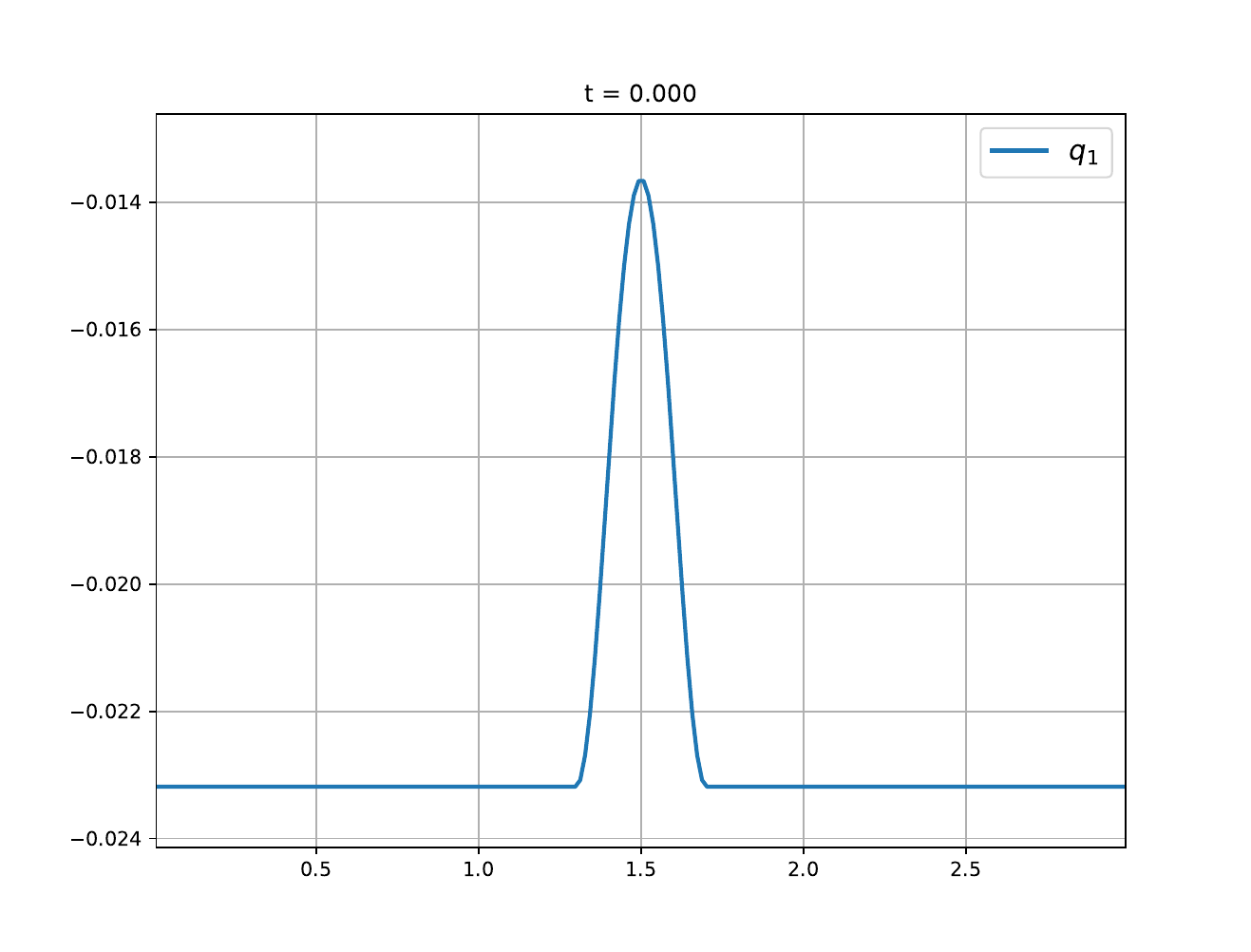}}
\subfloat[ $hu_2$ ]{
\includegraphics[width=0.5\textwidth]{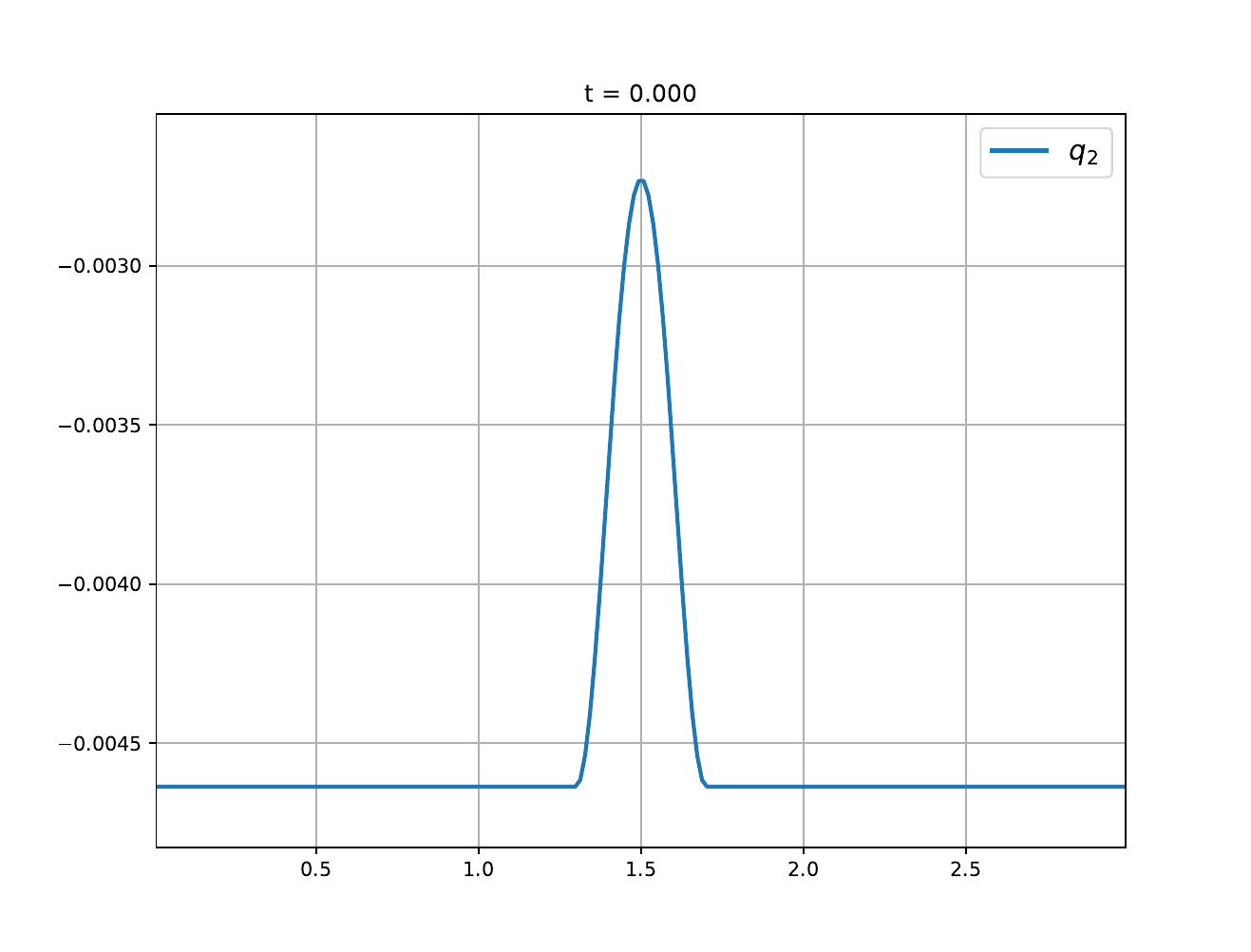}}
\end{center}
    \caption{Test 5. Perturbation of a subcritical steady state. Initial condition: free surface $\eta$ (up-left), $hu_0$ (up-right), $hu_1$ (down-left) and $hu_2$ (down-right).}
    \label{fig_test5_ci_perturb}
\end{figure}

We want to perform a test, similar to the one proposed in Test 4, 
but now considering a more physical initial condition. With that aim, we consider now $N=2$ moments, so that we have a parabolic velocity profile. The spatial domain is $[0,3]$, the bottom topography is again \eqref{test2_bottom} and a final time $t=0.1$ is considered in order to observe the propagation of the perturbation before it has left the domain. As initial condition, we consider a perturbation of a subcritical stationary solution. That subcritical steady state $U^{0,*}$  is characterised by the following constants:
\begin{equation}\label{test5_cini}
    C_1=0.5, \quad C_2=21.15525, \quad C_3=-0.005, \quad C_4=-0.001, 
\end{equation}
and the initial perturbation is then given by
\begin{equation}\label{test5_ci}
    U^0(x)=U^{0,*}(x)+\left(10^{-4}e^{-200(x-2)^2},0, \dots, 0 \right),
\end{equation}
and shown in Figure \ref{fig_test5_ci_perturb}. 

\begin{figure}[H]
 \begin{center}
     \subfloat[ Free surface $\eta$]{
\includegraphics[width=0.5\textwidth]{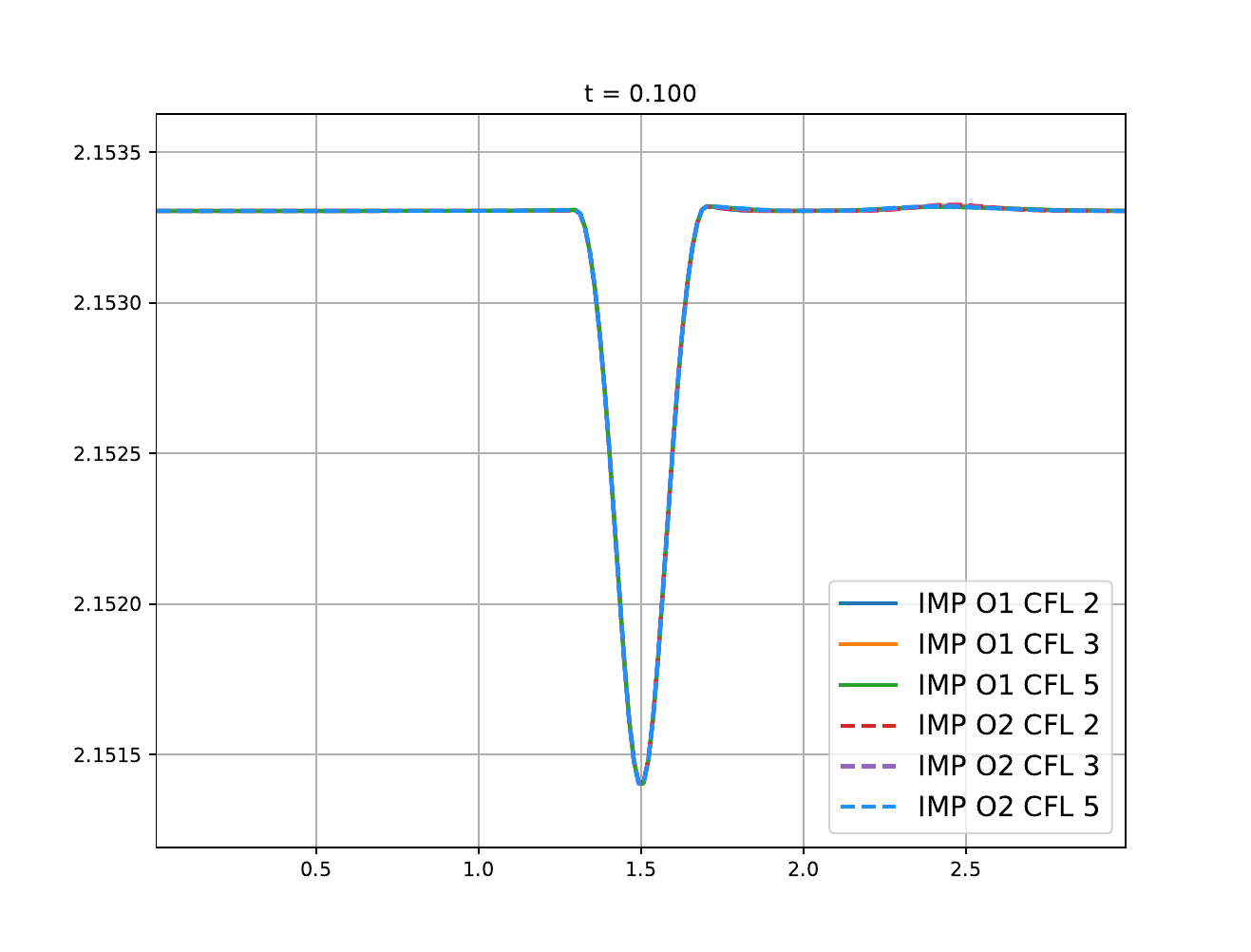}}
\subfloat[ Free surface $\eta$ (zoom) ]{
\includegraphics[width=0.5\textwidth]{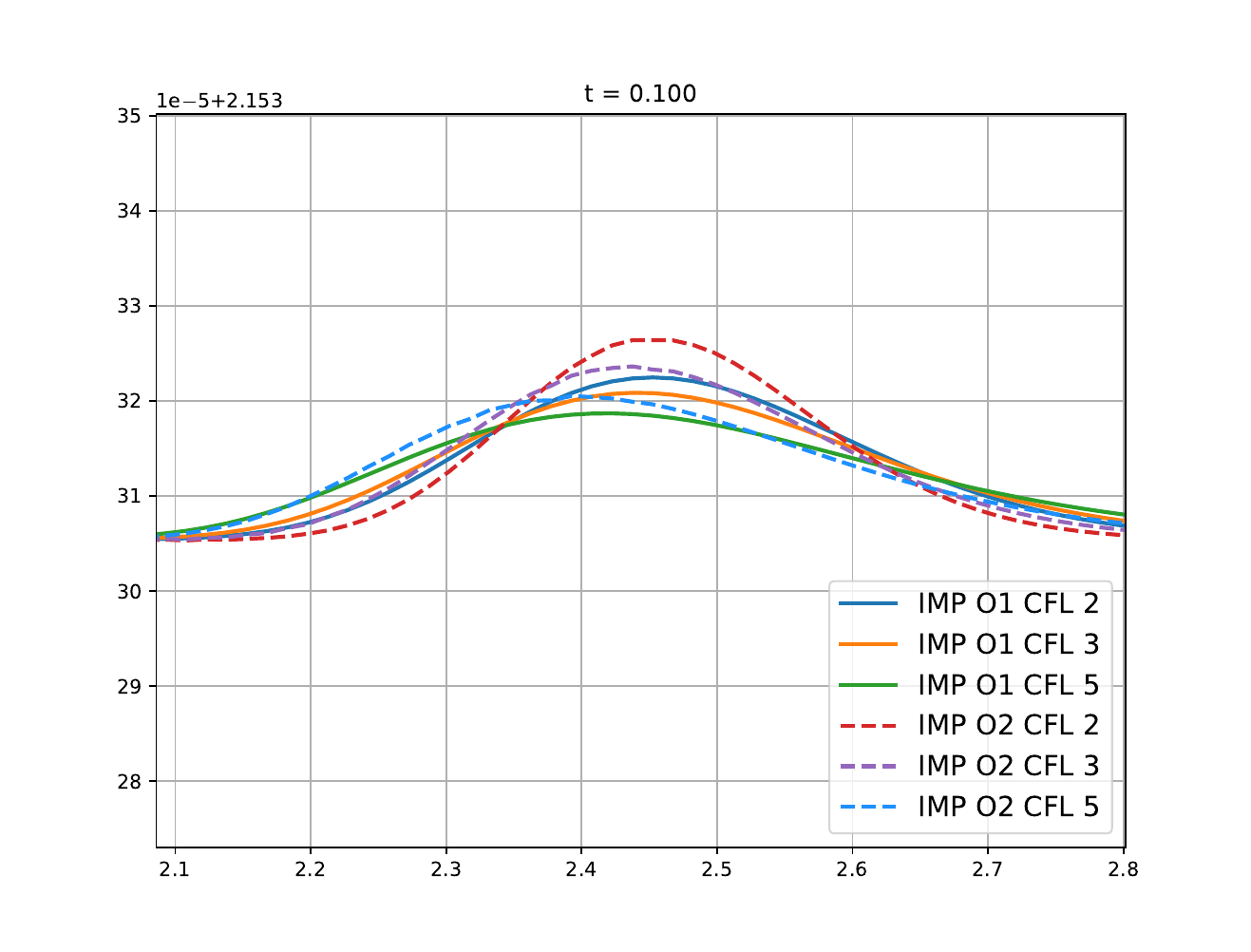}}
\vspace{2mm}
\subfloat[ $hu_0$ ]{
\includegraphics[width=0.5\textwidth]{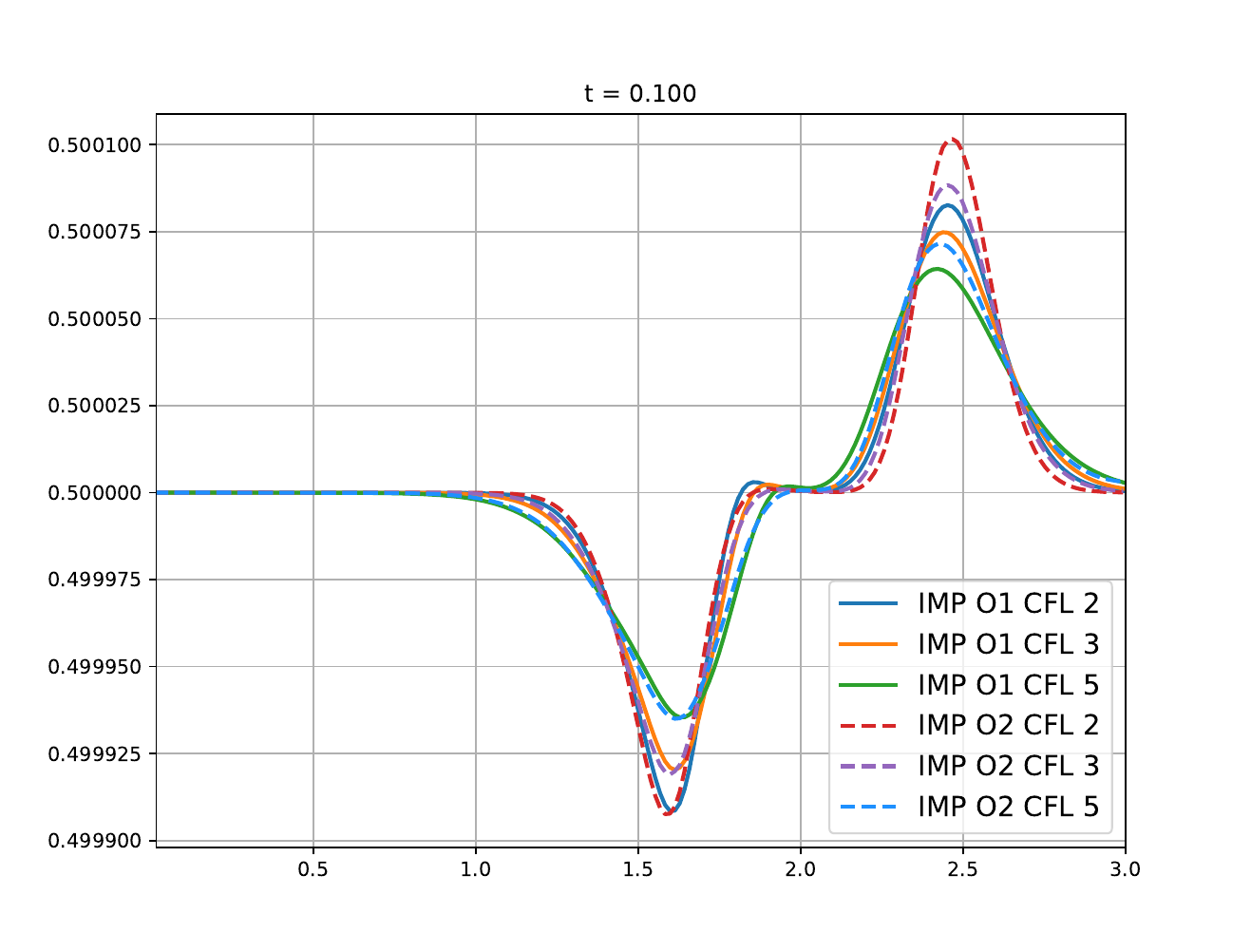}}
\subfloat[ $hu_1$ ]{
\includegraphics[width=0.5\textwidth]{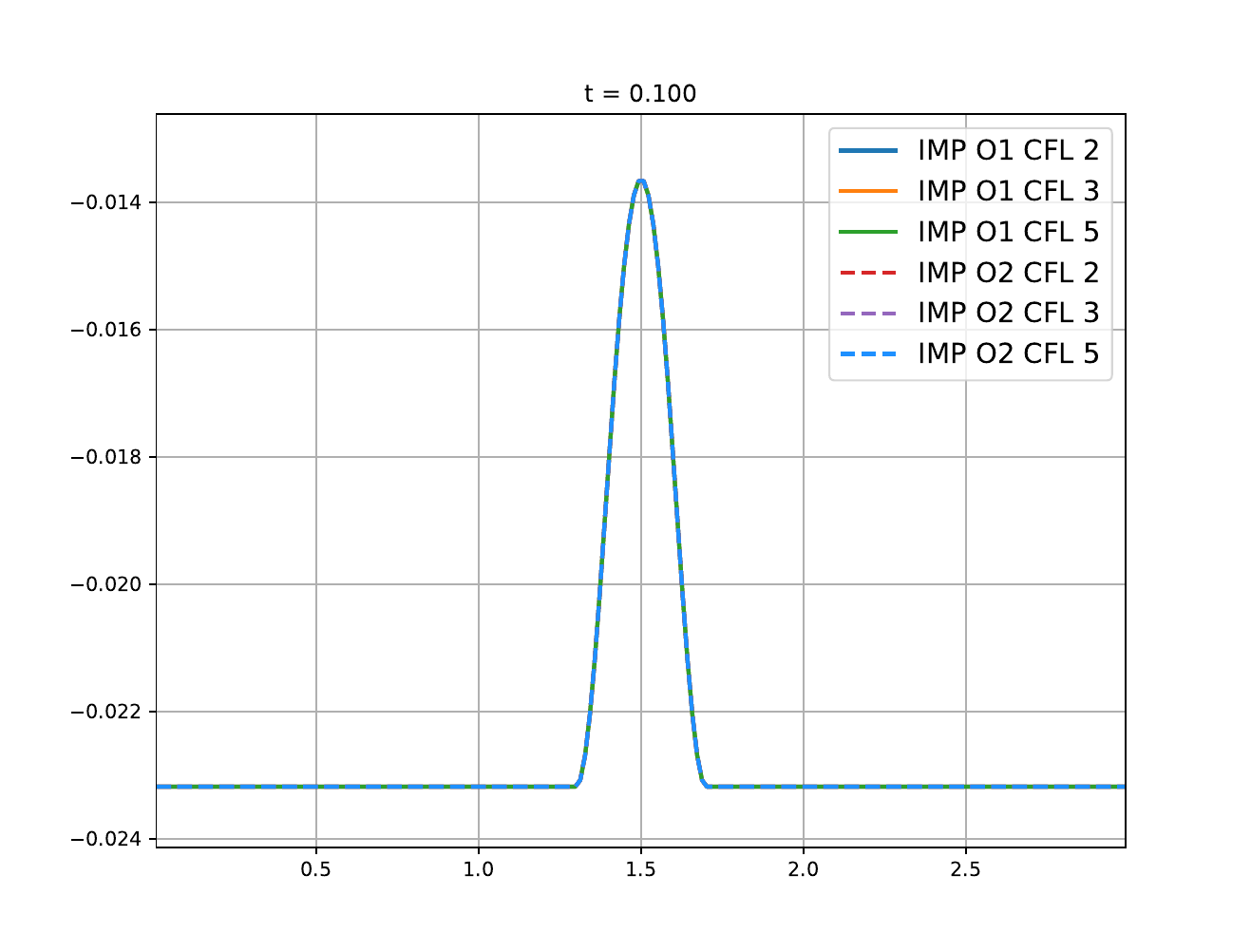}}
\vspace{2mm}
\subfloat[ $hu_2$ ]{
\includegraphics[width=0.5\textwidth]{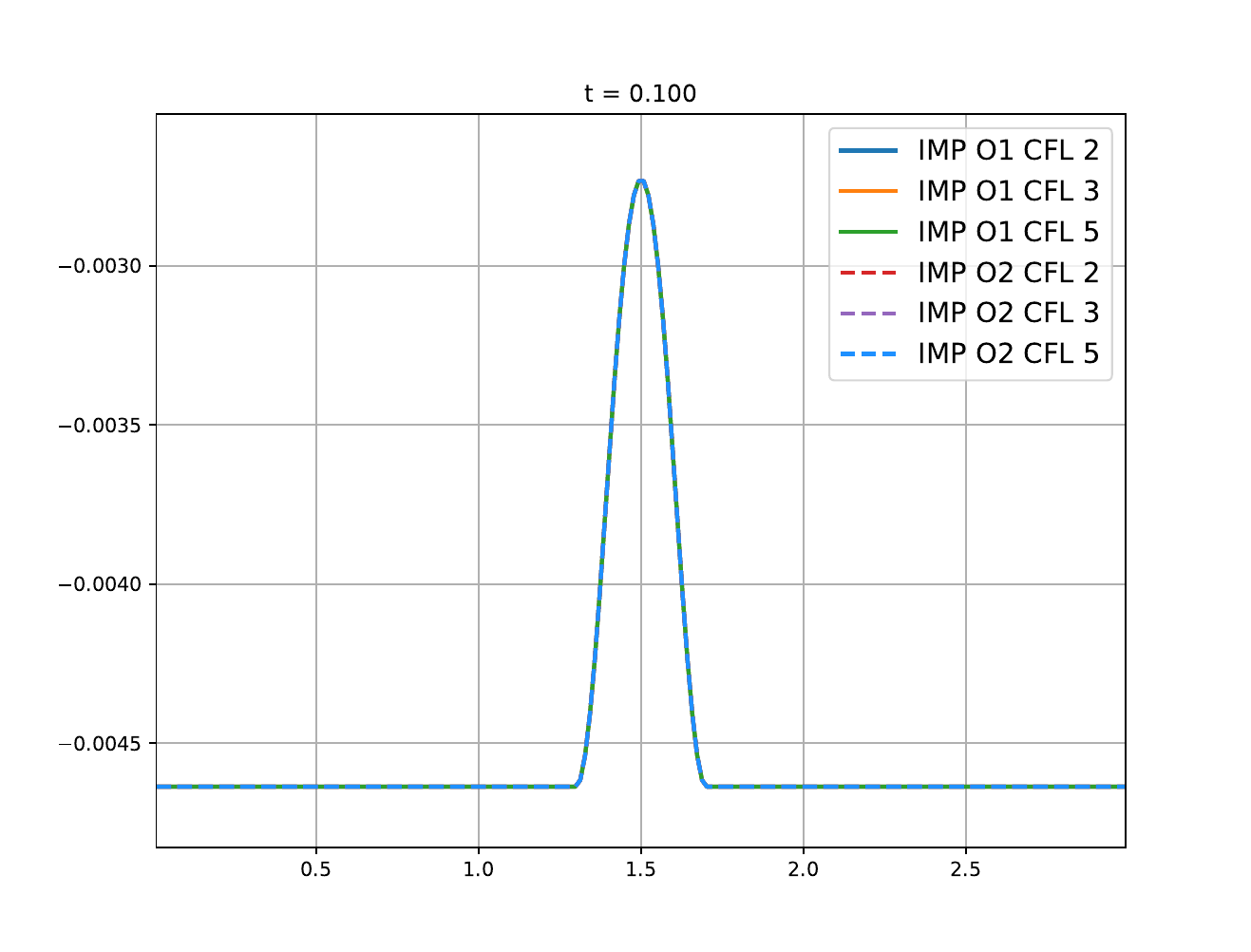}}
\end{center}
    \caption{Test 5. Perturbation of a subcritical steady state. Numerical solution for first and second order implicit schemes using different CFL values: free surface $\eta$ (up-left), zoom in the free surface (up-right), $hu_0$ (middle-left), $hu_1$ (middle-right) and $hu_2$ (down).}
    \label{fig_test5_perturb_t01}
\end{figure}

\begin{figure}[H]
    \centering
\includegraphics[width=0.75\linewidth]{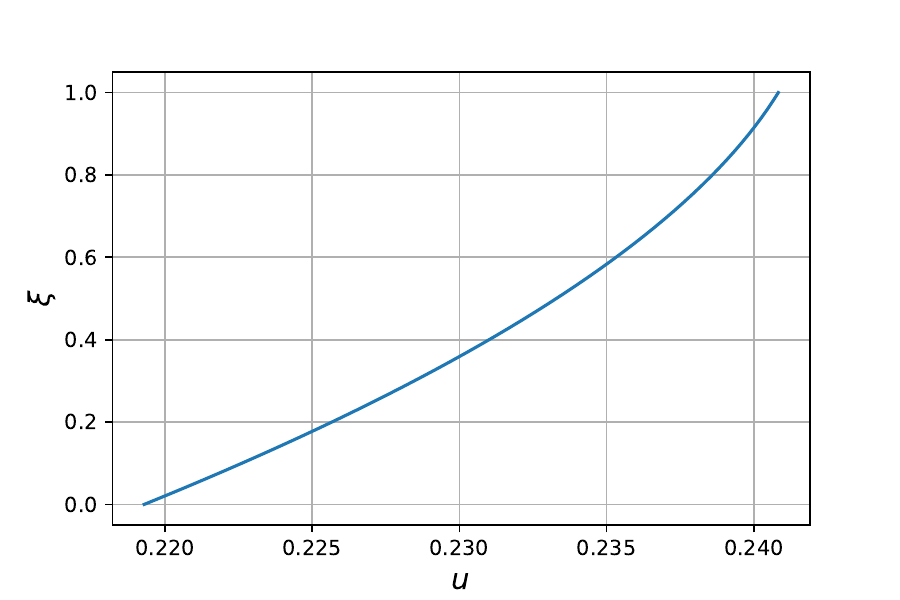}
    \caption{Test 5. Perturbation of a subcritical stationary flow: velocity profile at $x=2.25$ using the second order implicit scheme with CFL=5.} 
    \label{fig:test5_velocityprofile} 
\end{figure}

In Figure \ref{fig_test5_perturb_t01}, we have plotted the solution obtained with the implicit schemes considering different CFL values. We have zoomed in the free surface to see the differences between the schemes, observing, as expected that the first order schemes are, in general, more diffusive, and that the diffusion is increased when the CFL is bigger. This behaviour is also observed in the variable $hu_0$, while no major differences can be perceived in $hu_1$ and $hu_2$. As previously mentioned, the velocity profile is now parabolic. As an example, Figure \ref{fig:test5_velocityprofile} shows the vertical profile at the final time at point $x=2.25$ using the second order implicit scheme with CFL=5.

\subsection{Test 6. Dam-break problem}
Even though the aim of our schemes is to use them in low Froude number situations, since it is in that case when they outperform explicit schemes, we are now going to check that in fact they can also be used in other cases, such as in dam-break situations. In order to do so, we consider a problem similar to one shown in \cite{koellermeier2020analysis}. 
%However, since we are using the linearized model, in which higher-order moments are assumed to be small, those moments are considered here with smaller values than those considered for $u_0$. 
The bottom topography is considered to be flat, i.e., $z(x)=0$ and an initial condition given by: 
\begin{equation}\label{test6_cini}
\begin{split}
&h^0(x)=\left\{ \begin{array}{ll} 2 & \text{if }  x \leq 0, \\  
    1 & \text{otherwise,} \end{array} \right. \\
    & u_0^0(x)=0.25; \, u_1^0(x)=-0.005; \, u_8^0(x)=0.005; \, u_j^0(x)=0, \, j=2,\dots, 7.
    \end{split}
\end{equation}

The spatial domain  is $[-0.4,0.4]$ and the final time is $t=0.1$.
 
\begin{figure}[H]
 \begin{center}
\subfloat[ Free surface $\eta$]{
\includegraphics[width=0.5\textwidth]{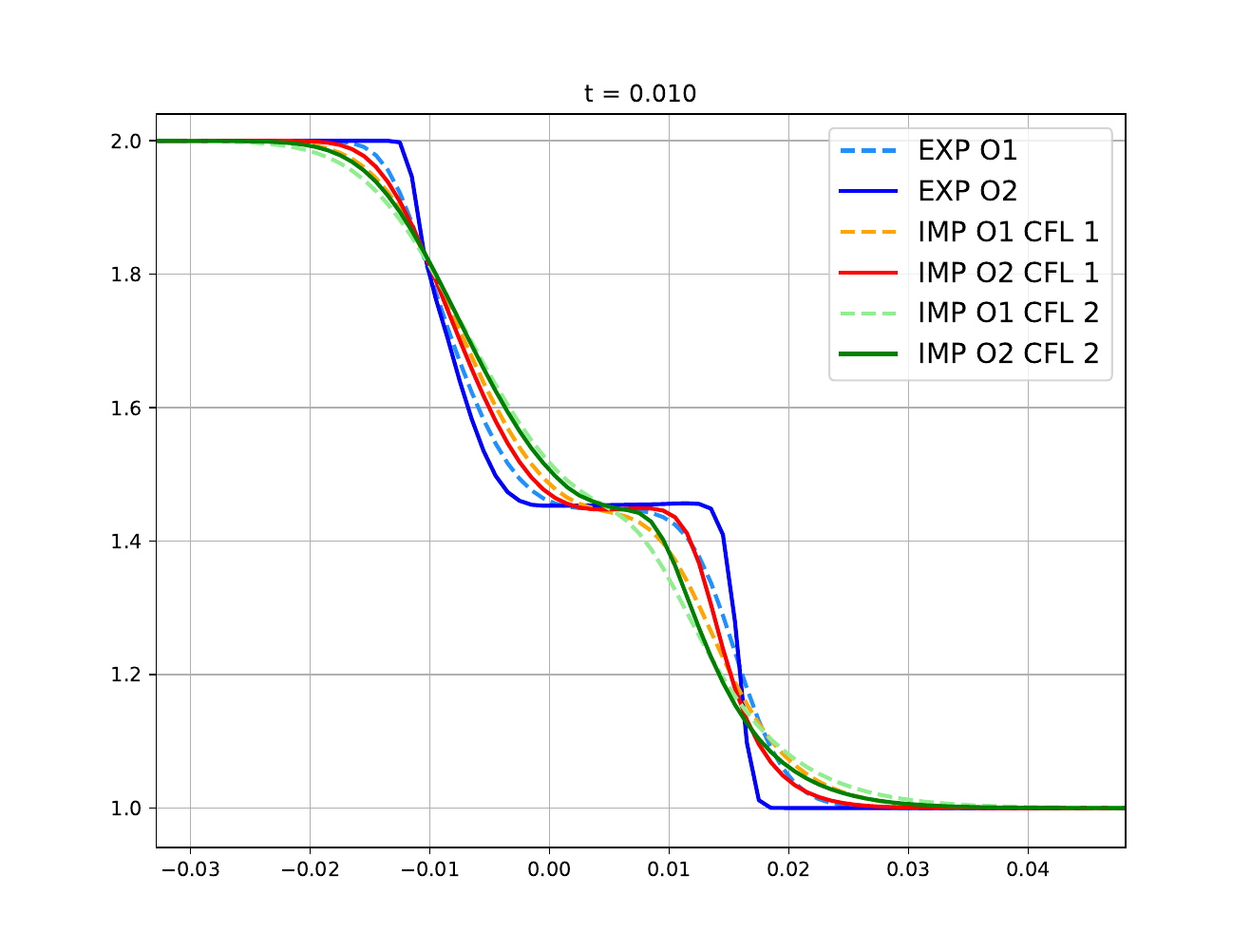}}
\subfloat[ Velocity $u_0$]{
\includegraphics[width=0.5\textwidth]{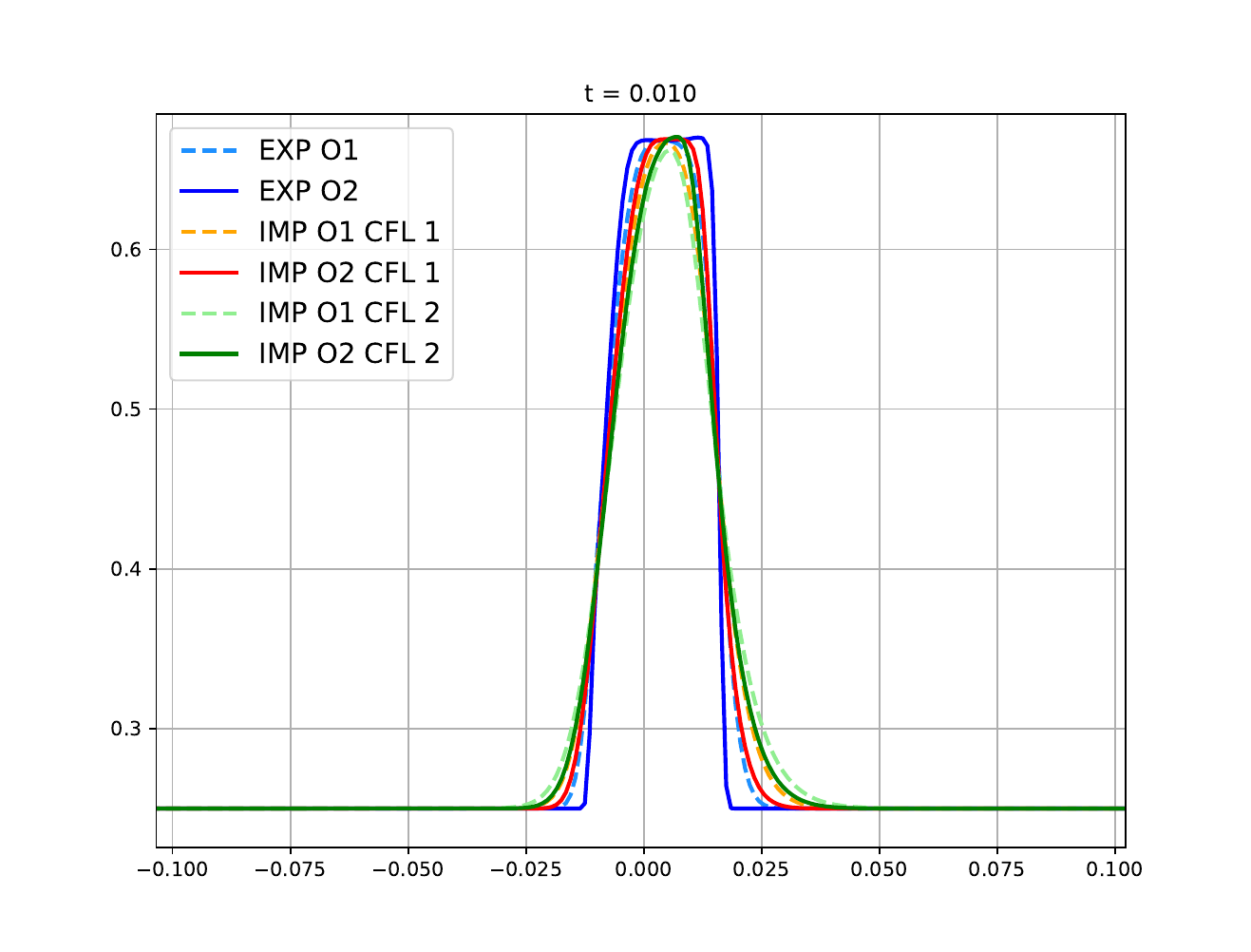}}
\vspace{2mm}
\subfloat[ First moment $u_1$]{
\includegraphics[width=0.5\textwidth]{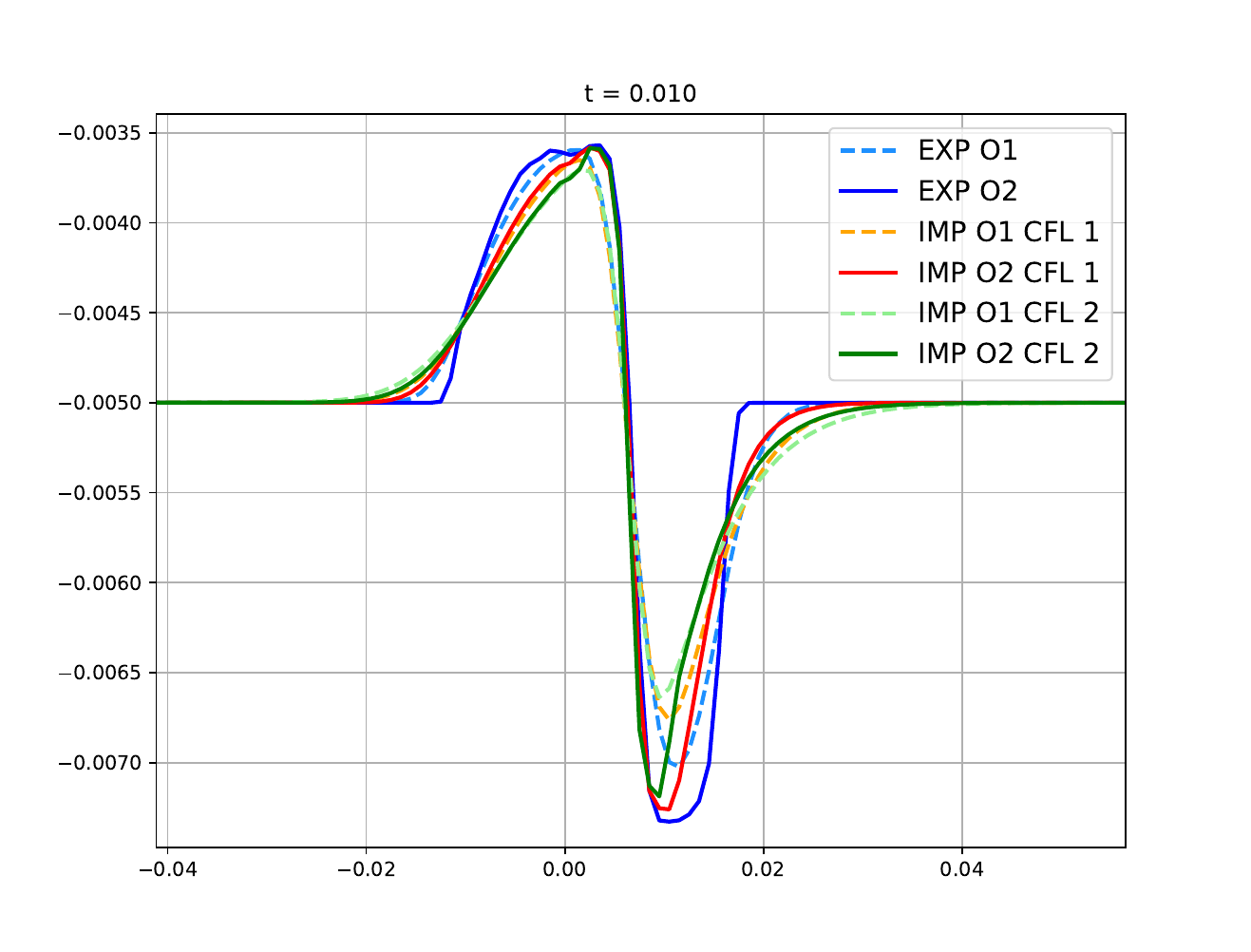}}
    \subfloat[ Last moment $u_8$]{
\includegraphics[width=0.5\textwidth]{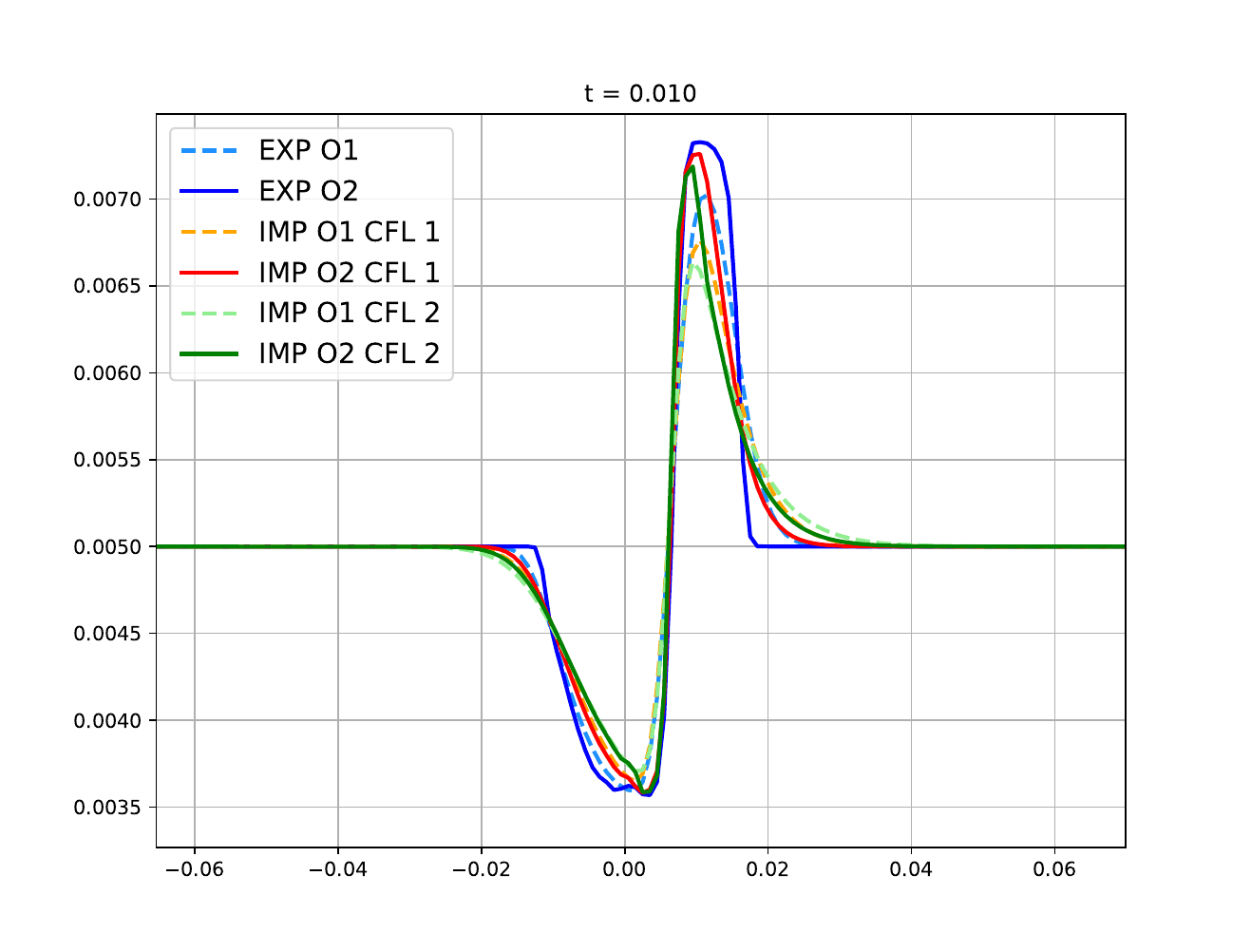}}
\end{center}
    \caption{Test 6. Standard dam-break test. Zoom of the variables $\eta$, $u_0$, $u_1$ and $u_8$ at time $t=0.01$. }
    \label{fig:test6_t001}
\end{figure}

\begin{figure}[H]
 \begin{center}
\subfloat[ Free surface $\eta$]{
\includegraphics[width=0.5\textwidth]{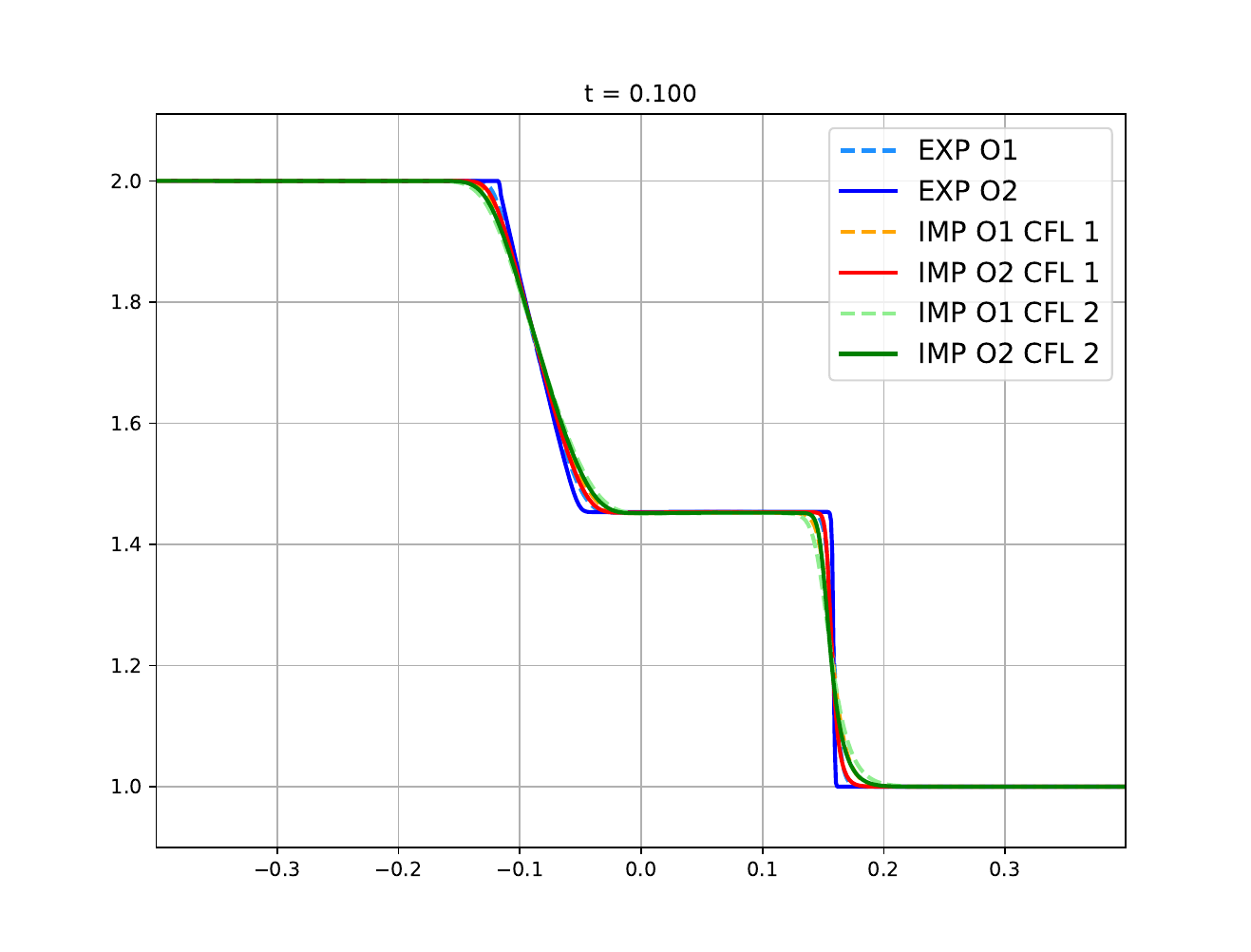}}
\subfloat[ Velocity $u_0$]{
\includegraphics[width=0.5\textwidth]{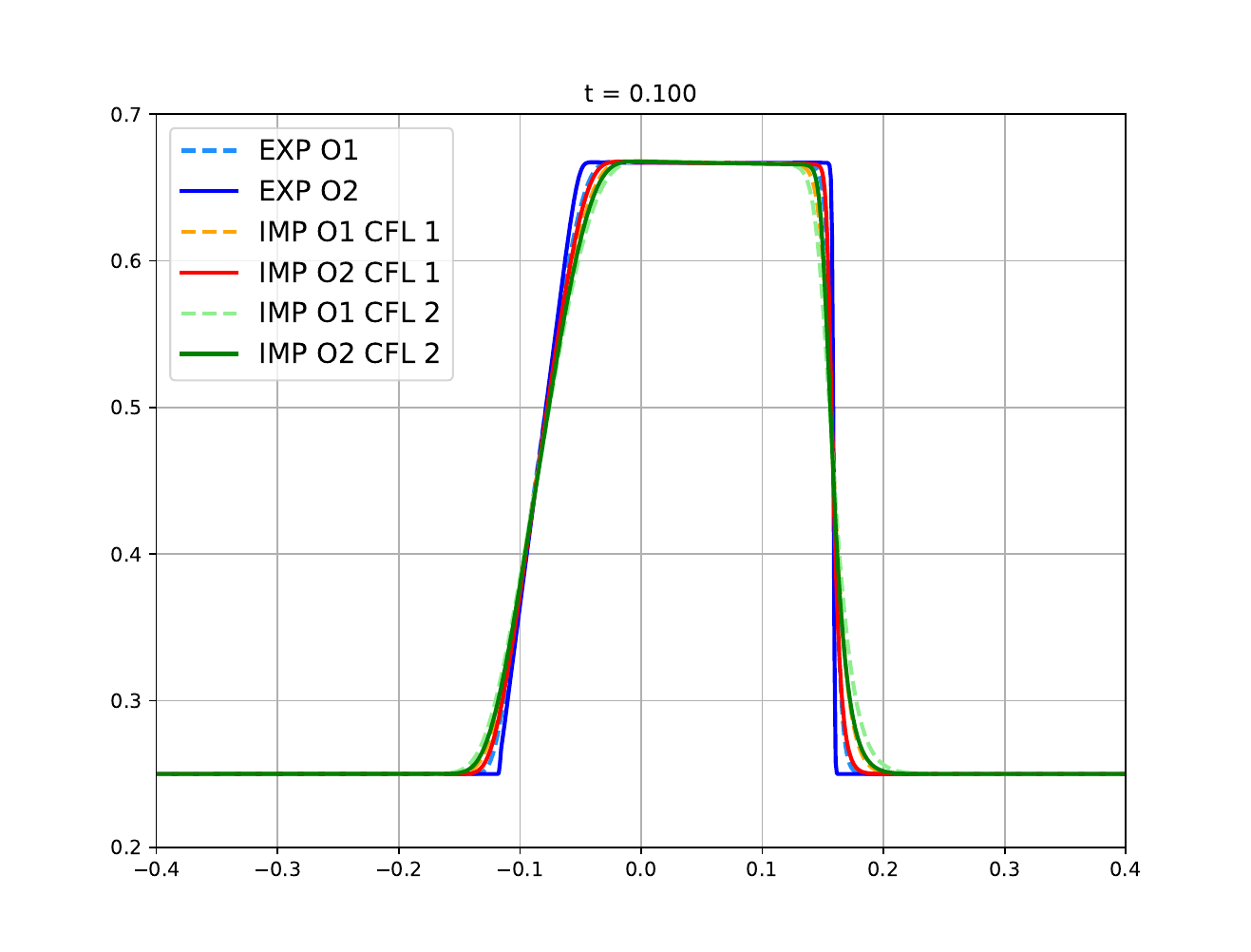}}
\vspace{2mm}
\subfloat[ First moment $u_1$]{
\includegraphics[width=0.5\textwidth]{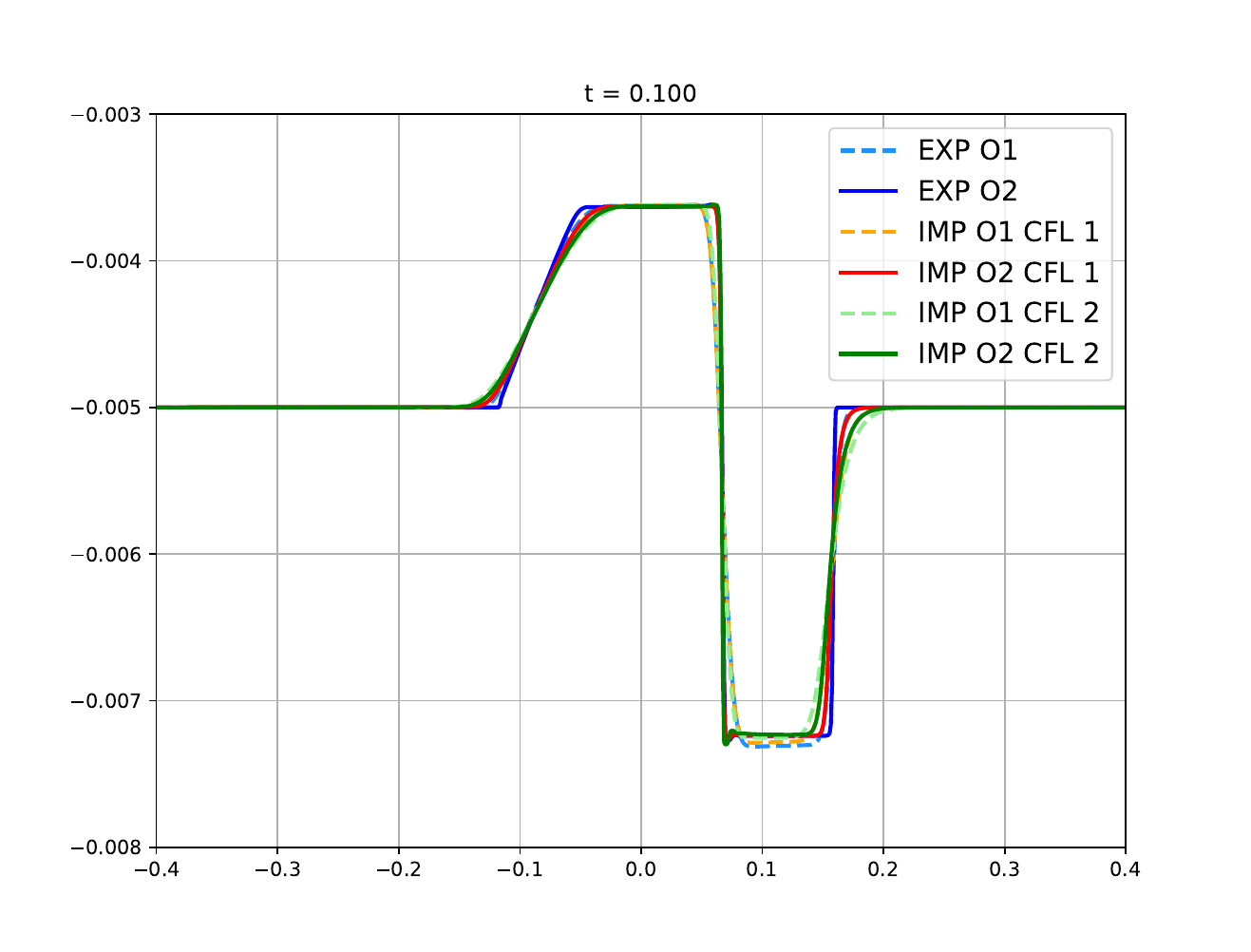}}
    \subfloat[ Last moment $u_8$]{
\includegraphics[width=0.5\textwidth]{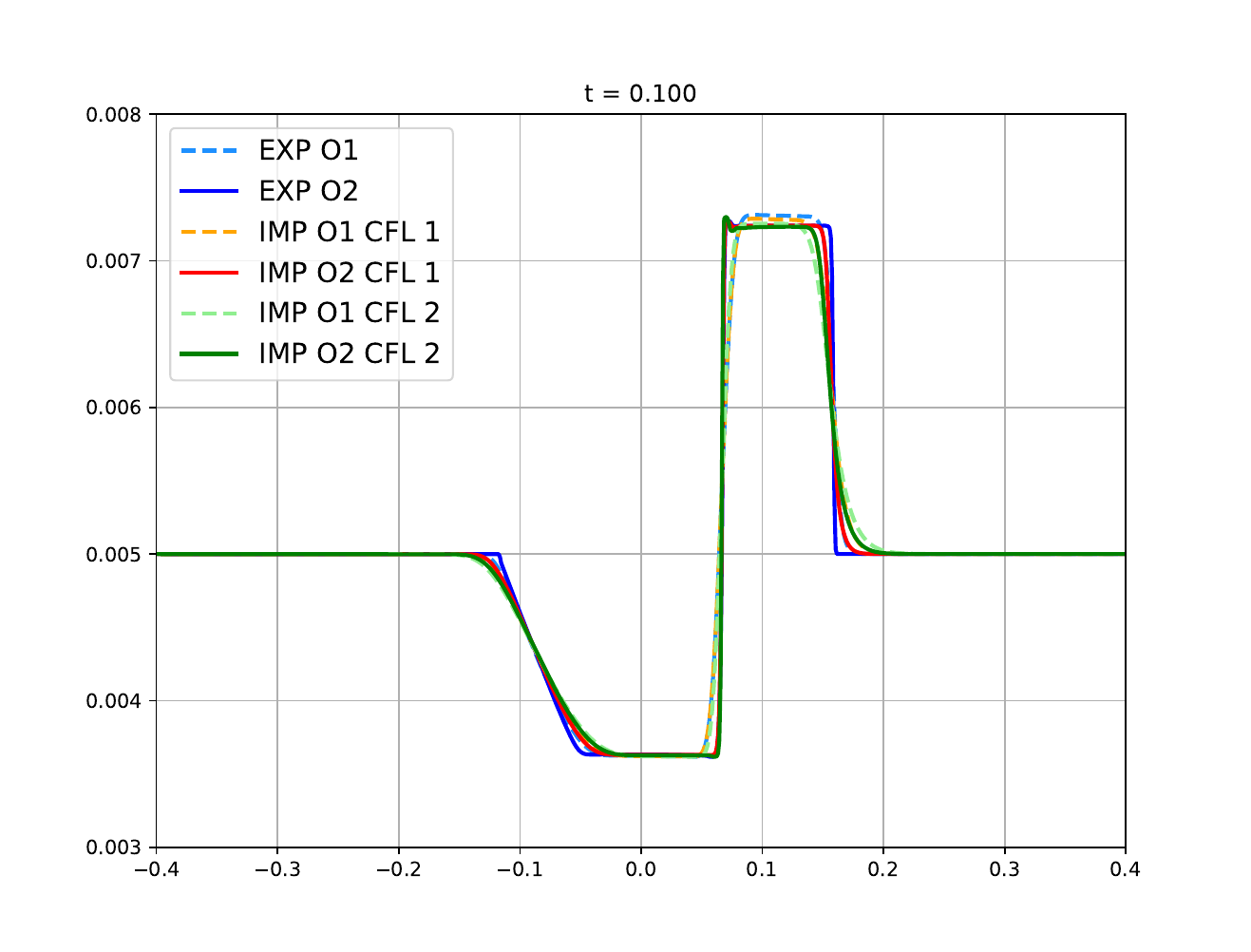}}
\end{center}
    \caption{Test 6. Standard dam-break test. Variables $\eta$, $u_0$, $u_1$ and $u_8$ at time $t=0.1$. }
    \label{fig:test6_t01}
\end{figure}

Figures \ref{fig:test6_t001} and \ref{fig:test6_t01} display the solutions provided by the methods of order 1 and 2 at $t=0.01$ and $t=0.1$ for $\eta$, the velocity $u_0$ and first and last moments $u_1$ and $u_8$. The CFL value is set to $0.9$ for the explicit schemes, whereas the values $1$ and $2$ (the second being the maximum allowed) are taken for the implicit schemes. 

\begin{figure}[H]
    \centering
    \subfloat[$x=0$]{
\includegraphics[width=0.5\linewidth]{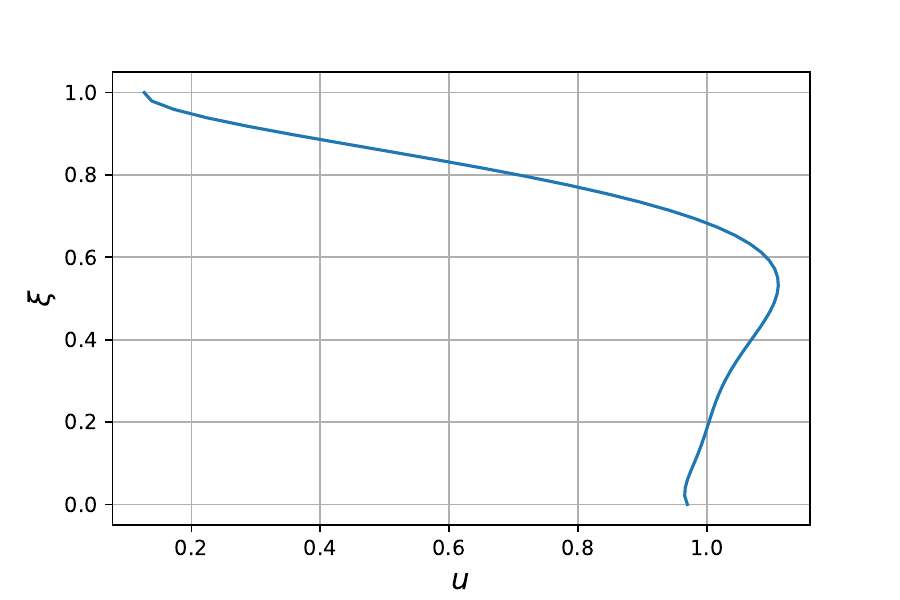}}
\subfloat[$x=0.15$]{
\includegraphics[width=0.5\linewidth]{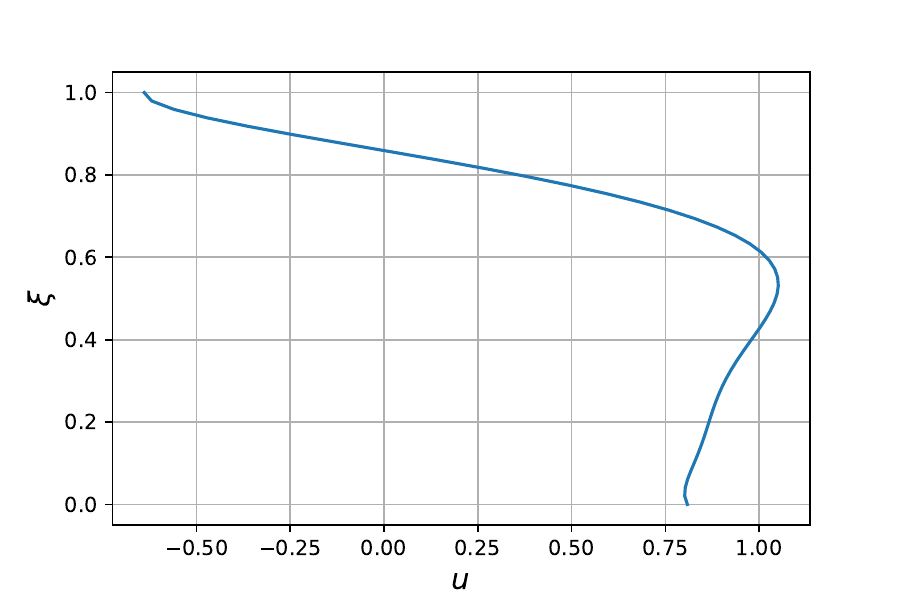}}
% \vspace{2mm}
% \subfloat[$x=0.35$]{
% \includegraphics[width=0.5\linewidth]{images/perfil_dambreak_tfinal_impo2_x035.pdf}}
    \caption{Test 6. Standard dam-break test: velocity profile at $x=0$ and $x=0.15$ at the final time $t=0.1$ using the second order implicit scheme with CFL$=2$. 
    %\textcolor{blue}{(JK: I wonder, are the profiles different for the other schemes? Since the values for $u_0, u_1, u_8$ are very similar, maybe it is almost the same. We decided not to plot the different profiles and wait for the reviewers. Maybe we can mention this as potential future works in the conclusion.)}
    } 
    \label{fig:test6_velprof}
\end{figure}

The images that depict moments can be reinterpreted as a vertical profile for velocity. As an example, in Figure \ref{fig:test6_velprof} we present the profiles obtained at the points $x = 0$ and $x = 0.15$ for the second order implicit scheme with a CFL number equal to 2 at the final time.

\section{Concluding remarks} \label{sec:Conclusion}
%JK: I like to repeat the abbreviations here for people who do not read the full paper and directly jump to the conclusion. But feel free to change it.
This work considers the development of numerical schemes for the Shallow Water Linearized Moment Equations (SWLME), which allow to obtain a vertical distribution of the horizontal velocity. This is in contrast to the standard Shallow Water Equations (SWE) where the horizontal velocity is constant in the vertical direction. The improved description of the velocity profile is useful in many practical applications, where a constant vertical profile could be very restrictive. Compared to other Shallow Water Moment Equations (SWME), the SWLME assume that higher order moments are small, which results in an explicit algebraic expression for steady states and explicit computation of the eigenvalues. 

SWLME have been solved by means of a semi-implicit first and second order exactly fully well-balanced method. A splitting approach allowed to separate acoustic and material phenomena. This is especially interesting for the case of low Froude regimes, where the CFL restriction due to the pressure term is too restrictive when compared to the velocity stability restriction. Indeed, using an implicit approach for the pressure system made it possible to use a larger CFL condition and a significant speed-up comparable to this increase is obtained. The resulting scheme is therefore efficient in such situations and exactly preserve the stationary solutions, as shown in the several numerical experiments that have been performed.

\section*{Acknowledgments}

This work is supported by proyect PID2022-137637NB-C21 funded by 
MCIN/AEI/10.13039/501100011033/ and ERDF A way of making Europe, and proyect PDC2022-133663-C21 funded by MCIN/AEI/10.13039/501100011033 and European Union NextGenerationEU/PRTR.

% Los primeros 9 polinomios de Legendre en función de \( y \), tras realizar el cambio de variable \( x = -2y + 1 \):

% 1. \( P_0(x) = 1 \)
%    \[
%    P_0(y) = 1
%    \]

% 2. \( P_1(x) = x \)
%    \[
%    P_1(y) = -2y + 1
%    \]

% 3. \( P_2(x) = \frac{1}{2}(3x^2 - 1) \)
%    \[
%    P_2(y) = 6y^2 - 6y + 1
%    \]

% 4. \( P_3(x) = \frac{1}{2}(5x^3 - 3x) \)
%    \[
%    P_3(y) = -20y^3 + 30y^2 - 12y + 1
%    \]

% 5. \( P_4(x) = \frac{1}{8}(35x^4 - 30x^2 + 3) \)
%    \[
%    P_4(y) = 70y^4 - 140y^3 + 90y^2 - 20y + 1
%    \]

% 6. \( P_5(x) = \frac{1}{8}(63x^5 - 70x^3 + 15x) \)
%    \[
%    P_5(y) = -252y^5 + 630y^4 - 560y^3 + 210y^2 - 30y + 1
%    \]

% 7. \( P_6(x) = \frac{1}{16}(231x^6 - 315x^4 + 105x^2 - 5) \)
%    \[
%    P_6(y) = 924y^6 - 2772y^5 + 3465y^4 - 2310y^3 + 735y^2 - 84y + 1
%    \]

% 8. \( P_7(x) = \frac{1}{16}(429x^7 - 693x^5 + 315x^3 - 35x) \)
%    \[
%    P_7(y) = -3432y^7 + 12012y^6 - 18018y^5 + 13860y^4 - 5544y^3 + 1155y^2 - 56y + 1
%    \]

% 9. \( P_8(x) = \frac{1}{128}(6435x^8 - 12012x^6 + 6930x^4 - 1260x^2 + 35) \)
%    \[
%    P_8(y) = 12870y^8 - 51480y^7 + 90090y^6 - 85800y^5 + 45045y^4 - 12600y^3 + 1785y^2 - 72y + 1
%    \]

% Estos son los resultados simplificados en función de \( y \), usando el cambio de variable \( x = -2y + 1 \).

\printbibliography

\end{document}